\documentclass[12pt]{article}
\usepackage{amsmath,amsfonts,amssymb,amsthm,amscd}
\usepackage{a4wide}
\title{\it {Relative Manin-Mumford for semi-abelian surfaces}
}
\author{D. Bertrand, D. Masser,  A. Pillay, U. Zannier
\footnote{{\it Authors' addresses} : bertrand@math.jussieu.fr,   david.masser@unibas.ch, A.Pillay@leeds.ac.uk, u.zannier@sns.it}}
\newtheorem{Theorem}{Theorem}
\newtheorem{Proposition}
{Proposition}
{Definition}
\newtheorem{Lemma}
{Lemma}
\newtheorem{Corollary}
{Corollary}
\newcommand{\C}{\mathbb C}
\newcommand{\R}{\mathbb R}
\newcommand{\Q}{\mathbb Q}
\newcommand{\Z}{\mathbb Z}

\newcommand{\G}{\mathbb G}

\def\mapright#1{\smash{\mathop{\longrightarrow}\limits^{#1}}}
\date{July  2013 (ArXiv version) \footnote{ {\it  AMS Classification} : 14K15, 12H05, 14K20, 11J95. {\it Key words} : semi-abelian varieties;  Manin-Mumford, Andr\'e-Oort and Zilber-Pink conjectures; differential Galois theory; polynomial Pell equations.}}

\begin{document}

\maketitle
\noindent
{\it Abstract}.\,- We show that Ribet sections are the only obstruction to the validity of the relative Manin-Mumford conjecture for one dimensional families of semi-abelian surfaces. Applications include special cases of the  Zilber-Pink conjecture for curves in a mixed Shimura variety of dimension four, as well as the study of polynomial Pell equations with non-separable discriminants.

 \tableofcontents
 
\section{Introduction}

\subsection{The viewpoint of group schemes : relative Manin-Mumford}
     
Let $\Q^{alg}$ be the algebraic closure of $\Q$ in $\C$, let $S$ be an irreducible algebraic curve over $\Q^{alg}$, and let  $G/S$ be a semi-abelian scheme over $S$ of relative dimension 2 and toric rank 1. We write $G_{tor}$ for the union of all the torsion points of the various fibers of $G  \rightarrow S $. This set $G_{tor}$ is also the set of values at all points of $S$ of the various torsion-sections of the group scheme $G/ S$. Let further $s : S \rightarrow G$ be a  section of $G/S$. The image of $s$ is an irreducible algebraic curve $s(S) = W$ in $G$,  defined over $\Q^{alg}$. Pursuing the theme of ``unlikely intersections" and relative versions of the Manin-Mumford conjecture (see \cite{Za}), we here  study the following  question, where ``strict" means ``distinct from $G$". 

\medskip
\noindent
{\bf Question 1} : {\it assume that $W \cap G_{tor}$ is infinite (i.e. Zariski dense in $W$). Must  $W$ then  lie  in a strict subgroup  scheme of $G/S$ ?  }

\medskip
Let us review some of the results  known along this line :

\smallskip

 i) the analogous question has a positive answer when $G/S$ is replaced by an abelian   scheme of relative dimension 2 : see \cite{MZ},  \cite{Za}, Theorem III.12 and \S 5 for non-simple ones, and \cite{MZbis} for the general case. 

ii) Assume that {\it the scheme $G/S$ is isoconstant}, i.e. isomorphic, after a finite base extension, to a product  $G = G_0 \times S$ with $G_0/\Q^{alg}$. Then, the Zariski closure $W_0$ of  the projection of $W$ to $G_0$ is an algebraic curve (or a point)  meeting $G_{0, tor}$ Zariski-densely.  By Hindry's generalization of Raynaud's theorem on the Manin-Mumford conjecture, see \cite{Hi}, \S 5, Thm. 2, $W_0$   is  a torsion-coset of a strict algebraic subgroup $H_0$ of $G_0$, 	and $W$ lies in a translate of $H = H_0 \times S$ by a torsion section. So, a positive answer is known in this case.

\medskip
Therefore, the first new case occurs  when $G$ is a non  isoconstant extension over $S$ of  an isoconstant elliptic scheme  by $\G_m $, i.e. when {\it $G$ is a semi-constant extension} in the sense of \cite{BP}. But as in \cite{BP} (though for different reasons), this case turns out to be more delicate, and the question can then  have a {\it negative} answer. A counterexample is given in  \cite{BE}, and we shall refer to the corresponding sections $s_R$ as  Ribet sections  of $G/S$ and to their images $W_R = s_R(S)$ as  Ribet curves. For a formal definition, see the end of this subsection (and \cite{JR} for its initial version); a more concrete description  is given in Appendix I.   In this paper, we will prove that in all other cases, our question has a positive answer; in other words :

\medskip
\noindent
{\bf Main Theorem  }. {\it Let $E/S$ be an elliptic scheme over the curve $S/\Q^{alg}$, and let $G/S$ be an extension of $E/S$ by $\G_{m/S}$. Let further $s : S \rightarrow G$ be a section of $G/S$, with image $W = s(S)$. 
 
 \smallskip
\noindent
{\bf (A)} Assume that $W \cap G_{tor}$ is infinite. Then,

i) either $s$ is a Ribet section;

ii) or  $s$ factors through   a strict  subgroup  scheme  of $G/S$.

\smallskip
\noindent
{\bf (B)} More precisely, $W \cap G_{tor}$ is infinite if and only if $~s~$ is  a Ribet section, or a torsion section, or  a non isoconstant section of a strict subgroup scheme of $G/S$.  }

\medskip
\noindent
We point out that this statement is invariant under isogenies $G \rightarrow  G'$ of the ambient group scheme, and under finite base extensions $S' \rightarrow S$. Throughout the paper, we will allow ourselves, sometimes tacitly, to perform such isogenies and base extensions.

\medskip

We now rephrase Part (A) of the theorem according to the various types of extensions $G/S$ and elliptic schemes $E/S$   which can occur, and explain in each case the meaning of ``isoconstant" in its Part (B)\,; a more concrete discussion of this array of cases is given in \S 2. Concerning the type of $E/S$, we recall that the scheme $E/S$ is isoconstant if and only if the $j$-invariant of its various fibers is constant; performing a finite base extension, we will then assume that $E$ is equal to $E_0 \times S$ for some elliptic curve $E_0/\Q^{alg}$. As for the type of $G/S$,

- either  $G$ is isogenous as a group scheme over $S$  to a direct product $\G_m \times E$. We then say that {\it the extension $G/S$ is isotrivial} and perform this isogeny.  Since $W$ is flat over $S$,  Conclusion (ii)  then reads : $W$ lies in a translate of $E/S$ or of  $\G_{m/S} = \G_m \times S$ by a torsion section  of $G/S$ (and Conclusion (i) does not occur); in this case, the isoconstant sections 	of the strict subgroup schemes are the translates by torsion sections of the constant sections of $\G_{m/S}$, or of  the constant sections of  $E/S$ if $E = E_0 \times S$ is constant.

- or the extension $G/S$ is not isotrivial. Conclusion (ii) then reads : $W$ lies in a translate of  $\G_{m/S}$ by a torsion section  of $G/S$; in this case, the isoconstant sections 	are the translates by torsion sections of the constant sections of $\G_{m/S}$.

\smallskip

Now,  whether $G/S$ is or is not an isotrivial  extension,

- we automatically get Conclusion (ii) if either  the scheme $E/S$ is not isoconstant, or it  is isoconstant, but $E_0$ does not admit  complex multiplications, or if $E_0$ is CM, but $G$ is isoconstant (a case already covered by \cite{Hi}, of course).

-  in the remaining case where $E/S$ is isoconstant, with a CM elliptic curve $E_0/\Q^{alg}$, and $G$ is not isoconstant (hence in particular, not isotrivial), Ribet sections   of $G/S$ do exist, their images  $W$ do not lie in strict subgroup  schemes of $G/S$ but meet $G_{tor}$ infinitely often, while not all sections $s$ satisfying the hypotheses of the theorem are Ribet sections.
\medskip

In other words, both conclusions (i) and (ii) of the Main Theorem  occur  in this last case  and are then mutually exclusive. However,  there is a way of reconciling them, through the setting of Pink's extension of the Andr\'e-Oort and Zilber  conjectures to mixed Shimura varieties, which we turn to in the next subsection \S 1.2 (see, e.g., Corollary 1).

\medskip
Another type of application of the Main Theorem is given in Appendix II of the paper~: this concerns the solvability of Pell's equations over polynomial rings, and extends some of the results of \cite{MZbis} to the case of non-separable discriminants.

 \medskip
To conclude this introduction, here is the promised  definition of a {\it Ribet section} $s_R : S \rightarrow G$ (see \cite{UIPB}, \S 1 - a more concrete, but analytic, characterization is given in Appendix I~; several other definitions are discussed in \cite{B-E}). In the notations of the Main Theorem, let   $\hat E /S$ be the dual of $E/S$: the isomorphism class of the $\G_m$-torsor $G$ over $E$  is given by a section $q : S \rightarrow \hat E$. Let  further ${\cal P} \rightarrow E \times_S \hat E$  be the Poincar\'e biextension of $E$
  and $\hat E$ by $\G_m$~:  by   \cite{De}, 10.2.13,  a section $s : S \rightarrow G$ of $G/S$ lifting a section $p : S \rightarrow E$ of $E/S$ is entirely described  by a trivialisation of the $\G_m$-torsor $(p, q)^*{\cal P}$  over $S$. Assume now that $E = E_0 \times S$ isoconstant, and admits complex multiplications, and let $f : \hat E_0 \rightarrow E_0$ be a non-zero antisymmetric isogeny (i.e., identifying $E_0$ and $\hat E_0$, a purely imaginary complex multiplication) which for simplicity, we here assume to be divisible by 2. Then, $(f(q), q)^*{\cal P}$ is a trivial  torsor in  a canonical way, and the corresponding  trivialization  yields a well-defined section $s = s(f)$ of $G/S$ above $p = f(q)$. When $G/S$ is semiconstant (which means~: when $q$ is not constant), any section $s_R$ of $G/S$ a non-zero multiple of which is of the form $s(f)$ for some antisymmetric $f$ will be called a Ribet section of $G/S$. So, on such a semi-abelian scheme $G/S$,   there exists essentially only one Ribet section $s_R$ (more precisely, all are linearly dependent over $\Z$),  and  by \cite{BE} (see also \cite{B-E}, and Appendix I),  its image   $W_R = s_R(S)$ meets $G_{tor}$ infinitely often. It  follows from Hindry's theorem (see Corollary 3 of \S 4) that the latter property characterizes Ribet sections among those sections of $G/S$ which project to a section  of $E/S$ of the form $p = f(q)$.

\subsection{The viewpoint of mixed Shimura varieties :  Pink's conjecture}
The consequences of the Main Theorem described in this section are discussed in \cite{UIPB}, which we here summarize for the convenience of the reader.

\medskip

Let $X$ be a modular curve, parametrizing isomorphism classes of elliptic curves with some level structure, let $\cal E$ be the universal elliptic scheme over $X$, with dual $\hat{\cal E}$, and  let $\cal P$ be the Poincar\'e bi-extension of ${\cal E} \times_X \hat{\cal E}$ by $\G_m$. This is a mixed  Shimura variety  of dimension 4, which  parametrizes points $P$ on extensions $G$ of elliptic curves $E$  by $\G_m$. A point of ${\cal P}(\C)$ can be represented by a triple $(E, G, P)$ , and is called special if the attached Mumford-Tate group is abelian, which is equivalent to requiring that $E$ has complex multiplications, that $G$ is an isotrivial extension, and that $P$ is a torsion point on $G$. Denote by ${\cal P}_{sp}$ the set of special points of $\cal P$. Following \cite{Pk}, we further say that an irreducible subvariety of $\cal P$ is  special if  it is a component of the Hecke orbit of a  mixed Shimura subvariety of  ${\cal P}$. Given any irreducible subvariety $Z$ of $\cal P$, the intersection of all the special subvarieties of $\cal P$ containing $Z$ is called the special closure of $Z$.  The special subvarieties of $\cal P$ of dimension $0$ are the special points; the special curves of ${\cal P}$ are described below;  for the  full  list, see  \cite{UIPB}, section 3. 

\begin{Corollary}  Let $W/\Q^{alg}$ be an irreducible closed algebraic curve in ${\cal P}$. Assume that $W \cap {\cal P}_{sp}$ is infinite. Then,  $W$ is  a  special curve.
\end{Corollary}

\medskip
To prove  this corollary, we distinguish the various  cases   provided by the projection $\varpi : {\cal P} \rightarrow X$ and its canonical section (rigidification) $\sigma : X \rightarrow \cal P$, whose image $\sigma(X)$ is made up of points of the type $(E, \G_m \times E, 0) \in {\cal P}$ :

\smallskip
- either the restriction of $\varpi$ to $W$ is dominant :  the corollary then says that $W$ lies in the Hecke orbit of the curve $\sigma(X)$. Indeed, up to Hecke transforms,  $ \sigma(X)$  is the only one-dimensional (mixed, but actually pure) Shimura subvariety of $\cal P$ dominating $X$.

In this case where $\varpi_{|W}$ is dominant, the  corollary follows not from our Main Theorem, but  from   Andr\'e's theorem \cite{AC}, p. 12,  on the special points of the mixed Shimura variety $\cal E$ (see also \cite{Pi0}, Thm. 1.2).

\smallskip

- or $\varpi(W)$ is a point $x_0$ of $X$, necessarily of CM type. In particular, $W$ lies in the fiber ${\cal  P}_0$ of $\varpi$ above $x_0$. This fiber ${\cal P}_0$ is a 3-dimensional mixed  Shimura subvariety of $\cal P$, which can be identified with the Poincar\'e biextension of $E_0 \times \hat{E_0}$ by $\G_m$, where $E_0$ denotes an elliptic curve in the isomorphism class of $x_0$.  An analysis of   the generic Mumford-Tate group of ${\cal P}_0$ as in \cite{Be}, p. 52, shows that up to Hecke transforms, there are exactly four types of special curves in ${\cal P}_0$ : the fiber   $(\G_m)_{x_0} $ above (0, 0) of the projection  ${\cal P}_0 \rightarrow ({\cal E}Ê\times_X \hat{\cal E})_{x_0} = E_0 \times \hat{E_0}$  and   the images $\psi_B(B)$ of the elliptic curves $B \subset E_0 \times \hat{E_0}$ passing through $(0,0)$ such that the $\G_m$-torsor ${\cal P}_{0|B}$ is trivial,  under the corresponding (unique) trivialization $\psi_B : B \rightarrow {\cal P}_{0|B}$. There are, up to torsion translates, three types of such elliptic curves $B$\;:  the obvious ones    $E_0 \times 0$ and  $0 \times \hat {E_0}$ (whose images we  denote by $\psi(E_0 \times 0), \psi(0 \times \hat {E_0})$), {\it and}  the  graphs of antisymmetric isogenies from $E_0$ to $ \hat{E_0}$, in which case $\psi_B$ corresponds precisely to a Ribet section of the semi-abelian scheme ${\cal G}_0 /\hat E_0$ defined below.

\medskip
Corollary 1 now follows from the Main Theorem, by interpreting ${\cal P}_0/\hat E_0$ as the ``universal" extension ${\cal G}_0$ of $E_0$ by $\G_m$, viewed as a group scheme over the curve $S := \hat E_0$, so that ${\cal P}_{sp} \cap {\cal P}_0 \subset ({\cal G}_{0})_{tor}$.  More precisely, suppose that $W$ dominates $\hat E_0$ : then, it is the image of a multisection of ${\cal G}_0/\hat E_0$, and after a base extension, the theorem implies that $W$ is the Ribet curve $\psi_B(B)$, or that it lies in a torsion translate of $\G_{m/\hat E_0} = \G_m \times \hat E_0$,  where  a new application of the theorem (or more simply, of Hindry's) shows that it must coincide with a Hecke transform of $\G_m = (\G_m)_{x_0} $ or of  $\psi(0 \times \hat E_0)$. By biduality (i.e. reverting the roles of $\hat E_0$ and $E_0$), the same argument applies if $W$ dominates $E_0$. Finally, if $W$ projects to a point of $E_0 \times \hat E_0$, then, this point must be torsion, and $W$ lies in the Hecke orbit of $(\G_m)_{x_0} $.  

\bigskip

Although insufficient in the presence of Ribet curves, the argument devised by Pink to relate the Manin-Mumford and the Andr\'e-Oort settings often applies (see the proof of  Theorems 5.7 and  6.3 of \cite{Pk}, and a probable addition to \cite{B-E} for abelian schemes). In the present situation, one notes that given a point $(E, G, P)$ in ${\cal P}(\C)$, asking that it be special as in Corollary 1  gives 4 independent conditions, while merely asking that $P$ be torsion on $G$ as in the Main Theorem  gives 2 conditions. Now, unlikely intersections for a curve $W$ in $\cal P$ precisely means  studying   its intersection with the union of the special subvarieties of $\cal P$ of  codimension $\geq 2$ (i.e. of  dimension $\leq 2$), and according to Pink's  Conjecture 1.2 of \cite{Pk}, when this intersection is infinite, $W$ should lie in a  special subvariety of dimension $ < 4$, i.e. a proper one. Similarly, if $W$ lies in the fiber ${\cal P}_0$ of $\cal P$ above a CM point $x_0$ and meets infinitely many special curves of this 3-fold, then, it should lie in a special surface of the mixed Shimura variety ${\cal P}_0$. In these directions, our Main Theorem,  combined with \cite{AC}, p. 12, and with  the relative version of Raynaud's theorem obtained in \cite{MZ},  implies  :  

\begin{Corollary}   Let $W/\Q^{alg}$ be an irreducible curve in the mixed Shimura 4-fold $\cal P$, and let $\delta_W$ be the dimension of  the special closure of $W$.  

i) Suppose that $\delta_W =   4$; then, the intersection of $W$ with the union of all the special surfaces of $\cal P$ {\it dominating $X$} is finite;

ii) Suppose that  $\delta_W = 3$; then, the intersection of $W$ with the union of all the special {\it non-Ribet} curves of $\cal P$  is finite.

\end{Corollary}

The proof goes along the same lines as that of Corollary 1; see \cite{UIPB} for  more details, and for a discussion on the gap between these corollaries and the full statement of Pink's conjecture 1.2 of \cite{Pk}, which would give a positive answer to  :

\medskip
\noindent
{\bf Question 2} {\it Let $W/\Q^{alg}$ be an irreducible curve in the mixed Shimura 4-fold $\cal P$, and let $\delta_W$ be the dimension of  the special closure of $W$. Is  the intersection of $W$ with the union of all the special subvarieties of $\cal P$ of dimension $\leq \delta_W - 2 $ then finite ?} 
 
 \medskip
 \noindent
 The case where $\delta_W = 2$ is covered by Corollary 1 . The remaining cases would consist in disposing of the restrictions   ``dominating $X$" and ``non-Ribet"  in Corollary 2.  These problems are out of the scope of the present paper. 
 
 \medskip
Although the Shimura view-point will not be pursued further, the Poincar\'e biextension, which has already appeared in the  definition of Ribet sections, plays a   role in the proof of the Main Theorem. See   \S 5.3, Remark 3.(iii), \S 3.3, Footnote (3), and  \S 6, Case ($\bf SC 2$), where   $s$ is viewed as a section of $\cal P$, rather than of $G$. See also the  sentence concluding \S 4.  

 \subsection {Plan of the paper}

\smallskip
 
 \noindent
$\bullet$ In the first sections (\S\S \,  2, 3.1, 3.2), we give a set of notations,   present the overall strategy of the proof, borrowed from \cite{MZ} and based on the same set of preliminary lemmas : large Galois orbits, bounded heights, Pila-Wilkie upper bounds. The outcome is that for a certain real analytic surface $\cal S$ in $\R^4$ attached to the section $s \in G(S)$:

\smallskip
\centerline{\bf $ W  \cap G_{tor}$  infinite $\Rightarrow  {\cal S}$ contains a semi-algebraic curve.}

\smallskip
\noindent
$\bullet$ The program for completing the proof is sketched in \S 3.3, and can be summarized by the following two steps. Here $log_G(s)$ is a local logarithm of $s$, and $F_{pq}$ is a certain Picard-Vessiot extension of $K$ attached to $G$ and to the projection $p$ of $s$ to $E$. Then, under a natural  assumption on $p$ (see Proposition 4),

\smallskip
\centerline{\bf ${\cal S}$ contains a semi-algebraic curve $\Rightarrow log_G(s)$ is algebraic over $F_{pq}$.}

\smallskip
\noindent
The proof is thereby reduced to a statement of algebraic independence, which forms the content of our ``Main Lemma", cf. middle of \S 3.3.  Notice the similarity between the statements of this Main Lemma and the Main Theorem, making it apparent that up to translation by a constant section,

\smallskip
\noindent
\centerline{\bf $ log_G(s)$ is algebraic over $F_{pq} \Rightarrow  s$ is Ribet or factors}

\smallskip
\noindent
as is to be shown.  

\smallskip
\noindent
$\bullet$ As a warm up, we realize these two steps in Section 4 when $G$ is an isotrivial extension.  In Section 5, we go back to the general case $G/S$, prove the first step, and comment on the use of Picard-Vessiot extensions.  

\smallskip
\noindent
$\bullet$ Section 6 is devoted to the proof of the Main Lemma. As in \S 4, we appeal to results of Ax type \cite{Ax} to treat  isoconstant cases, and to Picard-Vessiot theory, in the style of Andr\'e's theorem \cite{An}, for the general case.

\smallskip
\noindent
$\bullet$ Finally, Appendix I gives a concrete description of the local logarithm $log_G(s)$), including an analytic presentation of the Ribet sections $s_R$, and a proof, independent from those of \cite{BE}, \cite{B-E}, that they do contradict RMM (relative Manin-Mumford), while Appendix II is devoted to an application to polynomial Pell's equations, following the method of \cite{MZbis}.

\section{Restatement of the Main Theorem}
\subsection{Introducing $q$ and $p$}
     
We first repeat the setting of the introduction, with the help of  the fundamental isomorphism $Ext_S(E, \G_m) \simeq \hat E$ to describe the various cases to be studied.

\smallskip

So, let $\Q^{alg}$ be the algebraic closure of $\Q$ in $\C$. Extensions of scalars from $\Q^{alg}$ to $\C$ will be denoted by a lower index $_\C$. Let $S$ be an irreducible algebraic curve over $\Q^{alg}$, whose generic point we denote by $\lambda$, and  let $K = \Q^{alg}(S) = \Q^{alg}(\lambda)$, so $K_\C := \C(S_\C) = \C(\lambda)$. We use the notation $ {\overline \lambda}$ for the closed points of $S_\C$, i.e. $  {\overline \lambda} \in S(\C)$. Let  $G/S$ be an semi-abelian scheme over $S$ of relative dimension 2 and toric rank 1. Making Question 1 more precise, we write $G_{tor}$ for the union of all the torsion points of the various fibers of $G_\C \rightarrow S_\C$, i.e.  $G_{tor} = \cup_{{ {\overline \lambda}} \in S(\C)} (G_{ {\overline \lambda}})_{tor}  \subset G(\C)$, where $G_{\overline \lambda}$ denotes the fiber of $G$ above ${\overline \lambda}$. This set $G_{tor}$ is also the set of values at all ${{\overline \lambda}} \in S(\C)$ of the various torsion-sections of $G_\C \rightarrow S_\C$. Let finally $s : S \rightarrow G$ be a  section of $G/S$ defined over $\Q^{alg}$, giving a point $s(\lambda) \in G_\lambda(K)$ at the generic point of $S$, and closed points $s({ {\overline \lambda}}) \in G_{{\overline \lambda}}(\C)$ at the ${ {\overline \lambda}}$'s in $S(\C)$. 

\smallskip

In the description which follow,  we may have to withdraw some points of $S$, or replace $S$ by a finite cover, but will still denote by $S$ the resulting curve. After a base extension, the group scheme $G/S$ can be presented in a unique way as an extension 
$$0 \mapright{} \G_{m/S} \mapright{} 
G \mapright{\pi} E \mapright{} 0 ~$$
by $\G_{m/S} = \G_m \times S$ of an elliptic scheme $E/S$ over $S$. We denote by $\pi : G \rightarrow E$ the corresponding $S$-morphism. 

The extension $G$ is parametrized by a section 
$$q \in \hat E(S)$$ 
of the dual elliptic scheme $\hat E/S$. We write 
$$p = \pi \circ s  \in E(S)$$
 for the projection to $E$ of the section $s$. Finally, we denote by $k \subset \Q^{alg}$ the number field over which $S$ is defined, and assume without loss of generality that the generic fibers of $G$, hence of $q$, and of $s$, hence of $p$, are defined over $k(S) = k(\lambda)$.

\medskip
Since the algebraic curve $W = s(S) \subset G$ is the image of a section, the minimal strict subgroup schemes $H$ of $G$ which may contain  $W$ are flat over $S$, with relative dimension 1. These can be described as follows :  

-- either $q$ has infinite order in $\hat E(S)$ : then, $G/S$ is a non isotrivial extension, and  $H \subsetneq G$ is contained in a finite union of torsion translates of $\G_{m/S}$; in particular, $\pi(H) $ is finite over $S$, and $H$ can contain $W  $ only if $p = \pi(s)$ is a torsion section of $E/S$.

-- or $q$ has finite order : in this case, $G$ is isogeneous to the direct product $\G_{m/S} \times_S E$. Since the answer to the question is invariant under isogenies, we can then assume that $G$ is this direct product, i.e. that $q = 0$. Strangely enough, this (easy) case of Question 1 does not seem to have been  written up yet.  We present it in \S 4. The answer (cf. Theorem 1 below)  is a corollary of Hindry's theorem when $E/S$ is isoconstant, since $G/S$ too is then isoconstant; in an apparently paradoxical  way, we will use it  to characterize the Ribet sections of semi-constant extensions, cf. Corollary 3. 

\medskip

So, from now on and apart from \S 4, we could  assume that  $q$ has infinite order in $\hat E(S)$, i.e. that $G$ is a non-isotrivial extension.   However, the only hypothesis we will need   in our general study of \S 3.3 and \S 5 concerns the section $p = \pi(s)$ of $E/S$, see next subsection.

\subsection{Isoconstant issues}

In general, given our curve  $S/\Q^{alg}$, we say that a scheme $X/S$ is isoconstant if there exists a finite cover $S' \rightarrow S$  and a scheme $X_0/\Q^{alg}$ such that   $X_{S'} = X \times_S S'$ is isomorphic over $S'$  to the constant scheme $X_0 \times_{\Q^{alg}} S'$. We then say that a section $x$ of $X/S$ is {\it constant} if after a base change $S'/S$ which makes $X$ constant, the section of $X_{S'}/S'$ which $x$ defines comes from the constant part $X_0$ of $X_{S'}$. This notion is indeed independent of the choice of $S'$ (and often does not even require the base change $S'/S$).  See  Footnote (4) of \S 4 below for further conventions in the isoconstant cases.

\medskip

In these conditions, the hypothesis just  announced about $p$ reads :

\smallskip
\centerline{\bf $p \in E(S)$ is not a torsion section, and is not   constant  if $E/S$ is isoconstant,}

\smallskip
\noindent
and will be abbreviated as ``the section $p$ is not torsion, nor constant";  in terms to be described in \S5.3 and \S 6, it is better expressed as :

\smallskip
\centerline{\bf $p \in E(S)$ does not lie in the Manin kernel $E^\sharp$ of $E$.}

\smallskip
 
The relation to the Main Theorem is as follows. If $p$ is a torsion section, then, a torsion translate of $s \in G(S)$ lies in $\G_m$, so $s$ satisfies Condition (ii) of the Main Theorem. And   if $p = p_0  $ is a constant (and not torsion) section, then, $p(S) = \{p_0\} \times S$ does not meet $E_{tor}$ at all, so $s(S) \cap G_{tor}$ is empty. In other words, the Main Theorem is trivial in each of theses cases.

\medskip

 A precision on  the expression ``{\it $G$ is semi-constant}" is now in order   : it appears only  if $E/S$ is isoconstant, and  means that there exists a  finite cover $S' \rightarrow S$  such that  on the one hand,  there exists an elliptic curve $E_0/\Q^{alg}$ such that the pull-back of $E/S$  to $S'$   is isomophic over $S'$ to $E_0 \times S'$, and  that on the other hand, the section $q' \in \hat E'(S')$  which $q$ defines is given by a section of $\hat E_0\times_{\Q^{alg}} S$ which does {\it not} come from $\hat E_0(\Q^{alg})$. 
Since the answer to our question is invariant under finite base extension of $S$, we will assume in this case that $E/S$ is already constant, i.e. $E = E_0 \times S$, and  that $q \in \hat E_0(S) \setminus \hat E_0(\Q^{alg})$; indeed, as just said above, the second condition is then valid for any $S'/S$. Notice that the condition  that $q$ be non-constant forces it to be of infinite order. Consequently, semi-constant extensions are automatically both non isoconstant and non-isotrivial. On the other hand,   if $q \in \hat E_0(\Q^{alg})$ is constant, we are in the purely constant case of  \cite{Hi} already discussed in the Introduction.

\medskip
\noindent
{\bf Remark 1}  ({\it traces and images}) : let ${\cal E}_0/\Q^{alg}$ denote  the constant part ($K/\Q^{alg}$-trace) of $E/S$. An innocuous base extension allows us to assume that if $E/S$ is isoconstant, then  ${\cal E}_0 \neq \{0\}$ (and we then set ${\cal E}_0 :=   E_0$), and $E/S$  is actually isomorphic to $E_0 \times S$.  Denote by  $G_0/\Q^{alg}$ the constant part ($K/\Q^{alg}$-trace)  of $G/S$, and by $G^0/\Q^{alg}$ the maximal constant quotient ($K/\Q^{alg}$-image) of $G/S$. Then :

- $G/S$ is an isotrivial extension if and only if $q$ is a torsion section of $\hat E/S$.  In this case, an isogeny allows us to assume  that $q = 0$, i.e. $G = \G_m \times E$. We then have $G_0 = G^0 = \G_m \times {\cal E}_0$, and $G/S$ is isoconstant if and only if $E/S$ is isoconstant. 

-  assume now that $G/S$ is not an isotrivial extension, and that $E/S$, hence $\hat E/S$, is  not isoconstant.  Then, $G_0 = \G_m$, while $G^0 = {\cal E}_0 = \{Ê0\}$. 

- finally , assume that $G/S$ is not an isotrivial extension,  but $E = E_0 \times S$, hence $\hat E = \hat E_0 \times S$, is (iso)constant. Then either $q$ is a non constant section of   $\hat E/S$, in which case $G$ is semi-constant, and we have $G_0 = \G_m, G^0 =  {\cal E}_0 = E_0$; or $q$ is constant, in which case $G^0 = G_0 \in Ext_{\Q^{alg}}(E_0, \G_m)$ , and $G = G_0 \times S$ is itself constant.   
 
 \medskip
For the sake of brevity, we will henceforth say that a section of a group scheme $H/S$ is (iso)constant if it is a constant section of the constant part $H_0 \times S$ of $H/S$.  

\subsection{Antisymmetric relations and restatement}

The answer to Question 1 of the Introduction, as well as the proofs, depend  on possible relations between $p$ and $q$. For ease of notations,  we fix a principal polarization $\psi : \hat E \simeq E$ of the elliptic scheme, and allow ourselves to identify $q$ with its image $\psi(q) \in E(S)$. Also, we denote by $\cal O$ the ring of endomorphisms of $E$. If $E/S$ is not isoconstant, $\cal O$ reduces to $\Z$.  Otherwise, $\cal O$ may contain complex multiplications, and we say that the non-torsion and non-constant sections $p$ and $q$ are  {\it antisymmetrically related} if  there exists $\alpha \in {\cal O}\otimes \Q$ with  $\overline \alpha = - \alpha$ such that  $q = \alpha p$ 
in $E(K)$ modulo torsion.

\smallskip
Notice that we reserve this expression for sections $p, q$ which are non torsion (and non constant). Therefore,  an antisymmetric relation between $p$ and $q$  necessarily involves a non-zero imaginary $\alpha$,  hence ${\cal O} \neq \Z$,  forcing $E/S$ to be isoconstant.  And since $q$ is not constant,  the corresponding semi-abelian scheme $G= G_q$ is then semiconstant, and admits a Ribet section $s_R$ projecting to $p \in E(S)$.

\medskip

For any positive integer $m$, we set  :
$$S_m^G = \{{\overline \lambda} \in S(\C), s({\overline \lambda})~{ \rm has ~order}~  m ~{\rm in}~ÊG_{\overline \lambda} (\C)\} ,$$
$$S_m^E = \{{\overline \lambda} \in S(\C), p({\overline \lambda})~{ \rm has ~order}~ \  m ~{\rm in}~ÊE_{\overline \lambda} (\C)\} ,$$
$$S_{\infty}^G = \cup_{m \in \Z_{>0}} S_m^G \simeq W \cap G_{tor} ~, ~S_{\infty}^E = \cup_{m \in \Z_{>0}} S_m^E  \simeq  \pi(W) \cap E_{tor}   ,$$
where the indicated bijections are induced by $S \simeq s(S) = W \;,\; S \simeq p(S) = \pi(W) $. Clearly, $\cup_{k|m} S_k^G \subset \cup_{k|m}S_k^E$ for all $m$'s, and the points of $S_\infty^G$ can be described as those points of $\pi(W) \cap E_{tor}$ (``likely intersections") that lift to points of $W \cap G_{tor}$ (``unlikely intersections").

\medskip
Our Main Theorem can then be divided into the following three results. We first consider the case   when $q$ is torsion, which reduces after an isogeny to the case $q = 0$.

\begin{Theorem}    Let $E/S$ be an elliptic scheme over the curve $S/\Q^{alg}$, and let $G = \G_m \times E$ be the {\rm trivial} extension of $E/S$ by $\G_{m/S}$. Let further $s : S \rightarrow G$ be a section of $G/S$, with image $W = s(S)$,  such that $p = \pi(s)$ has {\rm infinite order in} $E(S)$. Then, $S_\infty^G$ is finite (in other words, $W \cap G_{tor}$ is finite) as soon as

(o) no multiple of $s$ by a positive integer factors through $E/S$ (i.e. the projection of $s$ to the $\G_m$-factor of $G$ is not a root of unity). 
\end{Theorem}

The case when $q$ is not torsion can be restated as follows. 

\begin{Theorem}  Let $E/S$ be an elliptic scheme over the curve $S/\Q^{alg}$, and let $G/S$ be a  non-isotrivial extension of $E/S$ by $\G_{m/S}$, i.e. parametrized by a section $q \in \hat E(S) \simeq E(S)$ {\rm of infinite order}. Let further $s : S \rightarrow G$ be a section of $G/S$, with image $W = s(S)$,  such that $p = \pi(s)$ {\rm has infinite order in} $E(S)$. Then, $S_\infty^G$ is finite (in other words, $W \cap G_{tor}$ is finite) in each of the following cases :

(i) $E/S$ is not isoconstant;

(ii) $E/S$ is isoconstant, and $p$ and $q$ are not antisymmetrically related;

(iii) $E/S$ is isoconstant,  $p$ and $q$ are  non-constant   antisymmetrically related sections,   and  no multiple of $s$ is a Ribet section.
\end{Theorem}

For the sake of symmetry, we recall that in these two theorems, the hypothesis that $p$ is not torsion is equivalent to requiring that no multiple of $s$ by a positive integer factors through $\G_{m/S}$. 

\medskip

Since   Ribet sections exist only in Case (iii)  of Theorem 2,
the conjunction of Theorems 1 and 2 is  equivalent to Part (A) of the Main Theorem, giving necessary conditions  for $W \cap G_{tor} \simeq S_{\infty}^G$ to be infinite. That these conditions  are (essentially) sufficient, i.e. that Part (B) holds true,  is dealt with by the following statement, which we prove right now.
\begin{Theorem} Let $E/S$ be an elliptic scheme over the curve $S/\Q^{alg}$, and let $G/S$ be an extension of $E/S$ by $\G_{m/S}$. Let further $s : S \rightarrow G$ be a section of $G/S$. Then,

i) if $s$ is a Ribet section,  $S_{\infty}^G$ is infinite, and equal to $S_{\infty}^E$;

ii) if $s$ is a torsion section, then $S_{\infty}^G  = S(\C)$;

iii)  if $s$ is a non-torsion section factoring through   a strict  subgroup  scheme  $H/S$ of $G/S$, then $S_{\infty}^G$ is empty if $s$ is an isoconstant section of $H/S$, and infinite but strictly contained in $S(\Q^{alg})$ if $s$ is not isoconstant.
\end{Theorem}
\noindent
{\it Proof} :  (i) is proved in \cite{BE}, \cite{B-E}, see also Appendix I below. The second statement is clear. As for (iii), this is an easy statement if the connected component   of $H$ is $\G_m$. If (for an isotrivial $G$), it is isogenous to $E$, the isoconstant case is again clear, while the non-isoconstant one follows from ``torsion values for a single point" as in \cite{Za}, cf. Proposition 1.(iv) below.

\medskip
So, we can now concentrate on Theorems 1 and 2.

\section{The   overall strategy}

Our stategy will be exacty the same as in \cite{MZ}. 

\subsection{Algebraic lower bounds}

\noindent
In this section, we denote by $k$ a number field over which the algebraic curve $S$, the group scheme $G/S$ and its section $s$ are defined. We fix a embedding of $S$ in a projective space over $k$, and denote by $H$ the corresponding  height on the set $S(\Q^{alg})$ of algebraic points of $S$.  We then have\;: 

\begin{Proposition} Let $E/S$ be an elliptic scheme, and let $p : S \rightarrow E$ be a section of $E/S$ {\rm of infinite order}. There exist  positive real numbers  $C  , C'$  depending only on $S/k, E/S$ and $p$, with the following properties.  Let ${\overline \lambda} \in S(\C)$ be   such that $p({\overline \lambda})$ is a torsion point of $E_{\overline \lambda}(\C)$, i.e. ${\overline \lambda} \in S^E_{\infty}$. Then,  

i) the point ${\overline \lambda}$ lies in $S(\Q^{alg})$, i.e. the field $k({\overline \lambda})$ is an algebraic extension of $k$;

ii) the height $H({\overline \lambda})$ of ${\overline \lambda}$ is bounded from above by $C$;

iii) if  $n \geq 1$ denotes  the order of $p({\overline \lambda})$, then $[k({\overline \lambda}):k]Ê \geq C' n^{1/3}$.

iv) the set $S^E_{\infty}$ is infinite (assuming that if $E/S$ is isoconstant, $p$ is not constant).
\end{Proposition}
\noindent
{\it Proof} : $i, ii, iii)$ In the  non-isoconstant case, one can reduce to the Legendre curve, where all is already written in \cite{MZ}, \cite{MZbis}, \cite{Za},  based on diophantine results of Silverman,  David    and the second author. Notice  that the  upper bound (ii) on $H(\overline \lambda)$  is needed to deduce the degree lower bound (iii). In the isoconstant case, the proof is easier as (ii) is not needed, and one can sharpen the lower bound (iii) to $n^2$ in the non-CM case (non effective), resp. $n/log n$ in the CM case. But as usual, any positive power of $n$ will do. 

$iv)$ This is the issue of ``torsion values for a single point", an analytic proof of which is given in \cite{Za}, p. 92, in the non-isoconstant case. If $E/S$ is isoconstant, the second case of this analytic proof does not occur, since we assume that $p$ is not constant.

\medskip
\noindent
{\bf Corollary}. {\it Hypotheses and the notation  $C'$ being  as in Proposition 1, let $G/S$ be an extension of $E/S$ by $\G_{m/S}$, and let $s$ be a section of $G/S$ lifting the non-torsion section $p$ of $E/S$.  Assume that ${\overline \lambda} \in S^E_{\infty}$ actually lies in $S^G_{\infty}$, and let $m$ be the order of the torsion point $s({\overline \lambda}) \in G_{\overline \lambda}(\Q^{alg})$. Then, $[k({\overline \lambda}):k]Ê \geq C' m^{1/6}$.}

\medskip
\noindent
{\it Proof} :  if $s({\overline \lambda})$ has precise order $m$ in $G_{\overline \lambda}$ and projects to a point $p({\overline \lambda})$ of order $n | m$ on $E_{\overline \lambda}$, then, $ns({\overline \lambda})$ is a primitive $m/n$-th root of unity, so has degree $>> (m/n)^{1- \varepsilon}$ over $k$, and we can assume that this is larger than $C'(m/n)^{1/3}$ (assume $C' < [k:\Q]^{-1}$). Since the fields of definition of $ns({\overline \lambda})$ and of $p({\overline \lambda})$ are contained  in $k({\overline \lambda})$, we get  $[k({\overline \lambda}):k] \geq C' (sup(n, m/n))^{1/3} \geq C'm^{1/6}$. 

\medskip
The conclusion of this first step is summarized by the implications :
$$ \forall {\overline \lambda} \in S(\C), {\overline \lambda} \in S^G_{\infty} \Rightarrow {\overline \lambda} \in S(\Q^{alg})Ê~{ \rm and }~ H({\overline \lambda}) \leq C,$$
and
$$ \forall {\overline \lambda} \in S(\Q^{alg}), \forall  m \geq 1,  {\overline \lambda} \in S_m^G \Rightarrow  [k({\overline \lambda}):k]Ê \geq C'm^{1/6}.$$
In particular, if $W = s(S)$ contains a point $w = s({\overline \lambda})$ of order $m$ (w.r.t. the group law of its fiber $G_{\overline \lambda}$), then $W$  contains at least $C'm^{1/6}$ points of order $m$ (w.r.t.  the group laws of their respective fibers) : indeed, since $W$ is defined over $k$, it  contains the orbit of $w$ under $Gal(\Q^{alg}/k)$. However, we will need a sharper version of  this statement, involving the archimedean sizes of the conjugates of $\overline \lambda$, and  the upper bound on $H({\overline \lambda})$ will again come to help at this stage.

\subsection{Transcendental upper bounds}
The next step is based on the following theorem of Pila. For the involved definitions and a short history on this type of results, leading to \cite{PW} and the accompanying references,  we refer to \cite{Za}, Remark 3.1.1,  and to the Amer. J. Math. version of \cite{MZ}). Dimensions here refer to real dimensions. For any $m \in \Z_{>0}$, we set $\Q_m = \frac{1}{m} \Z \subset \Q$.  

\begin{Proposition}
Let ${\cal S}$ be a naive-compact-2-(dimensional)-analytic subset of an affine space $\R^d$. Assume that no  semi-algebraic curve of $\R^d$ is contained in ${\cal S}$. For any $\varepsilon > 0$, there exists a real number $c =  c({\cal S}, \varepsilon) > 0$ with the following property. For each positive integer $m$, the set ${\cal S} \cap \Q^d_m$ contains  at most  $c\,m^\varepsilon$ points.
\end{Proposition}
\noindent
{\it Proof} : see \cite{MZ}, Lemma 1. 

\medskip
In \S 3.3, we will give a precise description of the real  surfaces $\cal S$ to which Proposition 2 is to be applied.  See the following \S 4.1  for a more easily recognizable form when  $G = \G_m \times E$.  Roughly speaking, 
$${\cal S}  \sim log^B_G(W) \subset (\Z \otimes  \C) \oplus (\Z^2 \otimes \R)  \simeq \R^4 $$
 is the set of logarithms of the various points $s({\overline \lambda})$, when ${\overline \lambda}$ runs through $S(\C)$, but we express these logarithms in terms of {\it a basis attached to the $\Z$-local system  of periods $\Pi_G$}, i.e. in Betti terms rather than in the de Rham view-point provided by the Lie algebra, hence the upper index $B$ above. In  this basis, the logarithms of the torsion points $s({\overline \lambda})$ of order $ m$  are represented by vectors with coordinates in $\Q_m$, and so $log^B_G(s(S_m^G)) \subset \Q_m^4$. Thanks to the ``zero estimate" discussed at the end of the present subsection, which compares the {\bf graph} $\tilde {\cal S}  \subset S \times  \R^4$ of $log^B_G \circ s : S \rightarrow  \R^4$  with its projection $\cal S$ in $\R^4$, Proposition 2 then implies  (with a proviso to be explained below for the first conclusion)
   
  \smallskip
 
-  either  that for any positive integer $m$,  $exp_G(\tilde {\cal S}) \sim s(S) = W$ contains at most  $c \;m^\varepsilon$ points of order $ m$  w.r.t the group law of their respective fibers;

-   or that $\cal S$ contains a real semi-algebraic curve, where algebraicity refers to the  real affine space $\R^4$ associated to   the above mentioned basis.

\smallskip
\noindent
In \S \S 4 and 5-6, we will prove that in all cases considered in Theorems 1 and 2,  the surface $\cal S$ contains no semi-algebraic curve. The first conclusion must then hold true. Combined with the conclusion of \S 3.1, this implies that the orders $m$ of the torsion points  lying on $W$ are uniformly bounded, and so, there exists a positive integer $N =  N(k,  S, G, s)$ such  that $W \cap G_{tor} \subset \cup_{{\overline \lambda} \in S(\Q^{alg})} G_{\overline \lambda}[N]$.  The latter set is  the union of the values at all ${\overline \lambda} \in S(\Q^{alg})$ of the torsion  sections  of $G/S$ of order dividing $N$, which form a finite union of   curves in $G$.  As soon as $p$ is not a torsion section, neither is $s$ and $W$ intersects this finite union in a finite number of points. Hence,  $W \cap G_{tor}$ is finite, and this concludes the proof of Theorems 1 and 2. (NB :   as done in \cite{MZ}, this conclusion can be  reached in a faster way  via the inequalities $m^{1/6} << [k(\overline \lambda) : k] <<    \sharp  (S_m^G) << m^\epsilon$, $H(\overline \lambda) << 1$, and the Northcott property.)

\medskip
\noindent
However, two points must be modified for the above discussion to hold.

\smallskip

$\bullet$ We need a uniform determination of the logarithms of  the points $s({\overline \lambda})$, and this requires fixing from the start a (any) simply connected pointed subset $(\Lambda, {\overline \lambda}_0)$ of the Riemann surface $S^{an}$ attached to $S(\C)$; in particular, our surface 
$${\cal S} = {\cal S}_\Lambda := log^B_{G, {\overline \lambda}_0}(s(\Lambda)),$$
and the graph $\tilde {\cal S}_\Lambda \subset \Lambda \times \R^4$ of $log^B_{G, {\overline \lambda}_0} \circ s : \Lambda \rightarrow  \R^4$ will depend  on a choice of $\Lambda$. Furthermore, the surfaces   $\cal S$ studied by Proposition 2   must  be compact,  so, $\Lambda$ too must  be compact. Consequently, $exp_G(\tilde {\cal S}_\Lambda) = s(\Lambda) \subset W(\C)$ is truly smaller than $W(\C)$, and the desired ``first conclusion" is reached in a slightly different way\;:  as in \cite{MZ} (see Lemma 6.2 of  AJM version, Lemma 8.2 of Math. Ann.), one first attaches to the height bound $C$ a finite union $\Lambda_C$ of simply-connected compact pointed sets  $(\Lambda_i, {\overline \lambda}_i)$ such that  for any point ${\overline \lambda} \in S(\Q^{alg})$ of sufficiently large degree with $H({\overline \lambda}) \leq C$, a positive proportion (say half, i.e. independent of ${\overline \lambda}$) of the conjugates of  ${\overline \lambda}$ over $k$  lies in $\Lambda_C$
\footnote{
~ In the present case, let us first remove from $S$ the finite set consisting of the points of bad reduction and those where the section is not defined (or any finitely many ones that may cause trouble along the way, possibly none...).   Now remove ``small" open disks around each of these points; what remains is a compact set in $S$. We want them small enough so that at most half of the conjugates of the relevant $\overline \lambda$  fall in their union: this may be achieved because $\overline \lambda$ has bounded height.
   In fact, if ``many" conjugates fall into a same small disk, then the corresponding contribution to the height is too big.
In turn this  follows for instance on looking at the difference $f(\overline \lambda) -f(\overline \lambda_0)$  where $\overline \lambda_0$ is the center of the disk and $f$  is  a suitable nonconstant coordinate on $S$.  Using the coordinate reduces the verification to the case of algebraic numbers  (rather than algebraic points). Having chosen these small enough disks, we cover the said compact set with finitely many  simply connected domains where the logs are locally defined.
}
.  Letting $\tilde {\cal S}_{\Lambda_C}$  be the finite union of the graphs  of  the maps $log^B_{G, {\overline \lambda}_i} \circ s$,  we deduce that for {\it any} ${\overline \lambda} \in S^G_m$, a similar positive proportion of the orbit of $s({\overline \lambda}) = w$ under $Gal(\Q^{alg}/k)$ lies in $exp_G(\tilde {\cal S}_{\Lambda_C}) = s(\Lambda_C) \subset W$. Proposition 2, combined with the subsequent zero estimate, then implies that this orbit has at most $c' m^\varepsilon$ points, and one can conclude as above (or  via the Northcott property on $\overline \lambda$). In what follows, we fix one of the  $\Lambda_i$'s, call it $\Lambda$ and write $log^B_G \circ s$ for $log^B_{G, {\overline \lambda}_i} \circ s$.

\smallskip

$\bullet$ Proposition 2 provides an upper bound for the image in ${\cal S}Ê\cap \Q_m^4 \subset \R^4$ of  $S^G_m \cap \Lambda$ under $log^B_G \circ s : \Lambda \rightarrow  \R^4$, while we need an upper bound for $S^G_m \cap \Lambda$ itself. In other words, we must show that not too many  points $\overline \lambda$ of $S_m^G$ can be sent by  $log^B_G \circ s$ onto the same point of $\Q_m^4$. Clearly, it suffices to show that  the projection
 $$u(\lambda) := d\pi(log_G(s(\lambda))) = log_E(\pi(s(\lambda))) =  log_E(p(\lambda))$$
 of $log_G(s)$ under the differential of $\pi : G \rightarrow E$ satisfies this separation property. The gap can now  be filled by appealing to the ``zero estimate" of \cite{MZ}, Lemma 7.1 of AJM, or Lemma 9.1 of Math. Ann.,   as follows.
 The (elliptic) Betti notations introduced here will be developed in \S 3.3 and \S 4.
 
 \smallskip
 
 Let $\Lambda$ be a compact  and simply connected subset of the Riemann surface $S(\C)$, which, without loss of generality, we may assume to be homeomorphic to a closed disk in $\C$. We recall that given an analytic sheaf $\cal F$ on $S(\C)$,  a section  $\sigma \in {\cal F}(\Lambda)$ of ${\cal F}$ over $\Lambda$ is by definition analytic on a neighbourhood of $\Lambda$ in $S(\C)$. Let now $E/S$ be an elliptic scheme over $S$,  and let $\omega_1(\lambda), \omega_2(\lambda)$ be the analytic functions on a neighbourhood of $\Lambda$ in $S(\C)$ expressing its periods relative to a given global  differential form of the first kind on $E/S$. Fix a determination $log_E$ of the corresponding elliptic logarithm on $E(\Lambda)$. For any analytic section ${\bf p} \in E(\Lambda)$, there then exists unique real analytic functions $\beta_1,\beta_2 : \Lambda \rightarrow \R$ such that $log_E( {\bf p} (\lambda)) = \beta_1(\lambda) \omega_1(\lambda) + \beta_2(\lambda) \omega_2(\lambda)$. We call $\{\beta_1, \beta_2\}$ the Betti coordinates of $\bf p$, and (extending a well-known notion for the Legendre family) say that $\bf p$ is a {\it Picard-Painlev\'e section} of $E/\Lambda$ if   its Betti coordinates $\beta_1, \beta_2$ are {\it constant}. We then have :

\begin{Proposition} {\bf (Zero estimate)}.  {\it Let $E/S$ and the simply-connected compact  subset $\Lambda$ be as above, and let $p$ be a regular section of $E/S$. Assume that   $p$ is not a torsion section, and if $E/S$ is isoconstant, that it is  not a constant section. There exists an integer $C''$ depending only on $E/S, \Lambda$  and $p$  such that for any Picard-Painlev\'e section $\bf p$  of $E/\Lambda$, the set
$$\{\overline \lambda \in \Lambda, log_E(p(\overline \lambda)) = log_E({\bf p}(\overline \lambda)\})$$
has at most $C''$ elements.}
\end{Proposition}

In other words, if we set  $u(\lambda) = log_E(p(\lambda)) = b_1(\lambda) \omega_1(\lambda) + b_2(\lambda) \omega_2(\lambda)$, then for any real numbers $\beta_1, \beta_2$, the equation $u(\overline \lambda) = \beta_1 \omega_1(\overline \lambda) + \beta_2 \omega_2 (\overline \lambda)$ has at most $C''$ solutions $\overline \lambda \in \Lambda$.  So, the statement above is just a fancy  translation of \cite {MZ}, Amer. JM, p. 12, and follows from Lemma 7.1 there in exactly the same way if $E/S$ is not isoconstant.  The constant case is even easier. Notice that we need Proposition 3 only for $ \beta_1, \beta_2$ running in $\Q $, i.e. for   torsion Picard-Painlev\'e sections $\bf p$, and that the Painlev\'e equation may bring a new view point on the computation of the bound $C''$. 

In our applications, $p$ is not a torsion section of $E/S$. And  in the isoconstant case, we have assumed without loss of generality that $p$ is not constant. So,   the map $log^B_E \circ p$ separates the points of $\Lambda$ up to the bounded error $C''$;  a fortiori, so does its lift $log^B_G \circ s$, and the gap between its image $\cal S$ and its graph  $\tilde {\cal S}$ is filled.

\subsection{What remains to be done}
In view of the previous discussion, the proof of   Theorems 1 and 2  is now reduced to defining the real surface $\cal S$ properly, and to showing that  under each of their hypotheses, $\cal S$  contains no semi-algebraic curve. This is dealt with as follows.

\subsubsection*{The real surface $\cal S$}

Fix a simply connected  and compact   subset $\Lambda \in S(\C)$,  homeomorphic to a closed disk, as well as a point ${\overline \lambda}_0$ in $\Lambda$, and a point $U_0$ in $ Lie(G_{{\overline \lambda}_0}(\C))$ such that $exp_{G_{{\overline \lambda}_0}}(U_0) = s({\overline \lambda}_0) \in G_{{\overline \lambda}_0}(\C)$. We henceforth denote by $\lambda$ the general element\footnote
{~ A   remark may be in order  about the meaning of the notation $\lambda$. In the first paragraphs, it represented the generic point of $S_\C$, i.e. we set $\C(S) = \C(\lambda)$ (by the way, from now on, we are over $\C$, so, dropping the lower index $\C$, we will write $K = \C(S)$). But it now represents the ``general element" of the simply connected domain $\Lambda \subset S^{an}$, which can have many analytic automorphisms.  It is understood that we here  consider only a global $\lambda$. Such a $\lambda$  may  require several  algebraically dependent parameters to be expressed.  Of course, when $\Lambda$ is small enough,  we can work with a chart $t$ on $\Lambda$ such that $\C(\lambda)$ is an algebraic extension of $\C(t)$. The results of functional  algebraic independence we appeal to do not require such reduction.}
 of $\Lambda$, and (sometimes) by an upper index $^{an}$ the analytic objects over the Riemann surface $S^{an}$ attached to our schemes over $S$. 
 
\bigskip
 We first repeat the definition of the real surface 
 ${\cal S}$ in more precise terms.
The group scheme $G/S$ defines an analytic family $G^{an}$ of Lie groups over the Riemann surface $S^{an}$. Similarly, its relative Lie algebra $(LieG)/S$ defines an analytic vector bundle $LieG^{an}$ over $S^{an}$, of rank 2. The $\Z$-local system  of periods of $G^{an}/\Lambda$  is the kernel of the  exponential exact sequence
$$0~ \mapright{} ~{\Pi}_G~ \mapright{} ~Lie G^{an} ~\mapright{exp_G}Ê~ G^{an} ~\mapright{}~ 0 ~, $$
 over $S^{an}$. Its sections over $\Lambda$ form a $\Z$-module $\Pi_G(\Lambda) \subset Lie G^{an}(\Lambda)$ of rank 3. Indeed, on using similar notations for the group schemes $E/S$ and $\G_m \times S$, 
the canonical projection $\pi : G\rightarrow E$ over $S$ induces at the Lie algebra level an exact sequence
$$0~ \mapright{} ~Lie \G_m^{an} ~ \mapright{} ~Lie G^{an} ~\mapright{d\pi}Ê~ Lie E^{an} ~\mapright{}~ 0 ~ .$$
From  the compatibility of the exponential morphisms, we deduce an exact sequence of $\Z$-local systems of periods
$$0~ \mapright{} ~{\Pi}_{\G_m}~ \mapright{} ~ \Pi_G  ~\mapright{d\pi}Ê~ \Pi_E  ~\mapright{}~ 0 ~, $$
 with ${\Pi}_{\G_m}(\Lambda)$ and $\Pi_E(\Lambda)$ of respective ranks 1 and 2 over $\Z$.
 
 \medskip
  
There exists a unique analytic section $U := \log_{G, {\overline \lambda}_0}$ of $Lie(G^{an})/\Lambda$ such that 
 $$U({\overline \lambda}_0) = U_0 ~Ê{\rm ~Êand~}Ê~\forall \lambda \in \Lambda, exp_{G^{an}_\lambda}(U(\lambda)) = s(\lambda). $$
 Since $\Lambda$ is fixed and ${\overline \lambda}_0$ plays no role in what follows, we will forget about them and will just write $U = log_{G, \Lambda} = log_G$, i.e. 
 $$\forall \lambda \in \Lambda, U(\lambda) = log_G(s(\lambda)).$$
 We call $U = log_G(s)$ ``the" logarithm of the section $s$. Its projection $p = \pi(s) \in E(S)$ admits as logarithm
 $$log_E(p ) := u  = d\pi(U) = d\pi( log_G(s)) = log_E(\pi(s)).$$ 
 We describe these logarithms in terms of classical Weierstrass functions in \S 4 for the (iso)trivial case $G = \G_m \times E$, and in Appendix I for the general case. These explicit expressions are not needed, but will provide  the interested reader with a translation of the algebraic independence results in more classical terms.
 
 \medskip
 Now, we  rewrite $U$ in terms of a conveniently chosen basis of  the $\Z$-local system of periods $\Pi_G/\Lambda$ of $G^{an}/\Lambda$. We call $U^B(\lambda) = log_G^B(s(\lambda))$ the resulting expression. For this, we choose a generator $\varpi_0 = 2\pi i$ of $\Pi_{\G_m}$, and a $\Z$-basis $\{\omega_1, \omega_2 \}$ of $\Pi_E(\Lambda)$. At each point $\overline \lambda \in \Lambda$, the latter generate over $\R$ the $\C$-vector space $Lie(E_{\overline \lambda})$. Consequently (and as already said before Proposition 3), there exist uniquely defined real analytic functions $b_1, b_2 : \Lambda \rightarrow \R^2$ such that 
 $$\forall \lambda \in \Lambda, ~u(\lambda) = b_1(\lambda) \omega_1(\lambda) + b_2(\lambda) \omega_2(\lambda). \qquad \qquad  ({\frak R}_u)$$
We call $u^B = (b_1, b_2) : \Lambda \rightarrow \R^2$ the Betti presentation of the logarithm $u = log_E(p)$.
  
 \medskip
 Now, choose at will lifts $\{\varpi_1, \varpi_2\}$ of $\{\omega_1, \omega_2\}$ in $\Pi_G(\Lambda)$. Then, $U - b_1 \varpi_1 - b_2 \varpi_2$ lies in  the kernel $Lie \G_m(\Lambda)$ of $d\pi$, which is generated over $\C$ by $\varpi_0$. Therefore, there exist a unique real analytic function $a : \Lambda \rightarrow \C = \R^2$ such that
 $$U = a \varpi_0 + b_1 \varpi_1 + b_2 \varpi_2.$$
 In conclusion, there exist uniquely defined  real-analytic functions
$a  : \Lambda \rightarrow \C, b_1 : \Lambda \rightarrow \R, b_2 : \Lambda \rightarrow \R$ such that
$U = log_G(s)$ satisfies the relation\;:
 $$\forall \lambda \in \Lambda, U(\lambda) = a(\lambda) \varpi_0(\lambda) + b_1(\lambda) \varpi_1(\lambda) + b_2(\lambda) \varpi_2(\lambda). \qquad ({\frak R}_U)$$
We call the  {\it real analytic}  map
$$ U^{B} =  (a,  b_1, b_2) : \Lambda \rightarrow \C \times \R^2 = \R^4,$$
the Betti presentation of the logarithm $U$ of $log_G(s)$. Its image
${\cal S} = {\cal S}_\Lambda := U^B(\Lambda) = log_G^B(s(\Lambda)) \subset \R^4$ is the real surface to be studied. Since $\Pi_G$ is the kernel of the exponential morphism, it is clear that for any $\overline \lambda \in S^G_m$, $U^B(\overline \lambda)$ lies in  $\Q_m   \times \Q_m^2 \subset \Q_m^4 \subset \R^4$.

\subsubsection*{Reducing to algebraic independence.}

To complete their proofs, we must show that under the hypotheses of Theorems 1 and 2, the surface $\cal S$ contains no semi-algebraic curve of the ambient affine space $\R^4$. This will be done in two steps, as follows. But before we describe them, we point out that since $log_E(p), log_E(s), \omega_i, \varpi_i, ...$ are local sections of the globally defined vector bundles $(Lie G)/S, (Lie E)/S$, it makes sense to speak of the minimal extension $K(log_G(s)), ...$ of $K = \C(\lambda)$ they generate  in the field of meromorphic functions over a neighbourhood $\Lambda'$ of $\Lambda$. A similar remark applies to the field $K(a), ...$ generated by the real analytic functions $a, ...$  in the fraction field of the ring of  real analytic functions over $\Lambda'$.

\medskip

$\bullet$ ({\it First step}) Let $F^{(1)}_{pq} = K(\omega_1, \omega_2, \log_E(p), log_E(q))$ be  the field   generated over $K$ by $\omega_1, \omega_2,  log_E(p) = u$ and a logarithm $v = log_E(q)$ of the section $q \in \hat E(S) \simeq E(S)$ parametrizing the extension $G$. As usual, we assume that if $E/S$ is isoconstant, then $p$ is not constant. We then claim that  {\it if  $\cal S$ contains a semi-algebraic curve,  and if $q = 0$, then,  $log_G(s)$ is algebraic over $F^{(1)}_{pq}$}.  In fact, we will only need a corollary of this result,  involving the universal vectorial extensions of $E$ and of $G$,  where  the base field $F^{(1)}_{pq}$ is replaced by a {\bf differential field} $F_{pq} := F^{(2)}_{pq}$ containing $F^{(1)}_{pq}$, which, inspired by the theory of one-motives, we can call the {\it field of  generalized periods of $\{E, p, q\}$}. (In more classical terms, the upper indexes (1) and (2) here stand for elliptic integrals of the first and second kinds.) In these conditions, and under {\it no} assumption on $q$, we will prove :

\begin{Proposition}  Assume that $p$ is not torsion and not constant, and that $\cal S$ contains a semi-algebraic curve. Then,    $log_G(s)$ is algebraic over the field  $F_{pq}$  of  generalized periods of $\{E, p, q\}$.
\end{Proposition}

\smallskip

$\bullet$ ({\it Second step}) The desired contradiction is then provided by the following Main Lemma, whose proof is the object of \S 6 (and \S 4 for $q=0$). This is a statement of ``Ax-Lindemann" type, but with logarithms replacing exponentials, in the style of Andr\'e's theorem \cite{An} (see also \cite{BM}) for abelian schemes. For results on semi-abelian surfaces close to this Main Lemma,  see \cite{Be}, Propositions 4.a, 4.b and Theorem 2.  For a broader perspective on algebraic independence of relative  periods, see Ayoub's recent paper \cite{Ay}.

\medskip
\noindent
{\bf Main Lemma}.  {\it With $S/\C$, let $G/S$ be an extension by $\G_m$ of an elliptic scheme $E/S$, parametrized by a section $q$ of $\hat E/S$, and  let $G_0$ be the constant part of $G$. Let further $s $ be a section of $G/S$, with projection $p = \pi \circ s $ to $E/S$,  and  let $F_{pq}$ be the  field of  generalized periods of $\{E, p, q\}$. 

\smallskip
\noindent
 {\bf (A)} Assume that   $log_G(s)$ is algebraic over $F_{pq}$. Then,  there exists a constant section $s_0 \in G_0(\C)$ such that

\smallskip

i) either  $s - s_0$ is a Ribet section;

ii) or  $s - s_0$ factors through a strict subgroup scheme of $G/S$.

\smallskip
\noindent
{\bf (B)} More precisely, $log_G(s)$ is algebraic over $F_{pq}$ if and only if  there exists a constant section $s_0 \in G_0(\C)$ such that  $s - s_0$ is a Ribet section, or a torsion section, or  factors through a strict subgroup scheme  of $G/S$ projecting onto $E/S$. }

\medskip

The analogy with the Main Theorem is clear, except perhaps for   the last conclusion of Part (B) of the Main Lemma (which forces an isotrivial $G \simeq \G_m \times E$). This is due to the fact that even in this isotrivial case,  the roles of $\G_m$ and $E$ are here not symmetric,  because of the occurence of $p$ in the base fields $F^{(1)}_{pq}, F_{pq}$. On the contrary (``torsion values for a single point" on a group scheme of relative dimension 1 over a curve), they  played similar roles for relative Manin-Mumford.

\medskip
\noindent
{\bf Proposition 4 $+$ Main Lemma (A) $\Rightarrow$ Theorems 1 and 2}.

\smallskip
Let us first deal with Theorem 1, where $G = \G_m \times E$,  with constant part $G_0 = \G_m$, resp. $G_0 \times S = G$, if $E$ is not, resp. is, isoconstant.   Assume for a contradiction that $S_\infty^G$ is infinite. Then, the real surface $\cal S$ must contain an algebraic curve, and since $G$ admits no Ribet section, Proposition 4, combined with Part (A) of the Main Lemma, implies that a multiple by a non-zero integer of the section $s$ factors through a translate of $H = \G_m$ or of $H = E$ by a constant  (non-necessarily torsion) section $s_0 \in G_0(\C)$. But the projection $p$ of $s$ to $E$ is by assumption not torsion, and we know that it cannot be constant.  So, $H$ must be equal to $E$, and $s$ projects on  the $\G_m$-factor of $G$ to a constant point $\delta_0$. Since $s(S) = W$ contains torsion points, $\delta_0$ must be a root of unity, and  $s$ factors through a torsion translate of $E$, as was to be proved.

\smallskip
In the direction of Theorem 2, we now assume that $G$ is a non-isotrivial extension, so $H = \G_{m/S}$ is the only connected strict subgroup scheme of $G/S$. The proof above  easily   adapts to  the case when $G$ is isoconstant, where again, $G$ admits no Ribet section (in the sense of \S 1). Now, assume that $G$ is not isoconstant, so $G_0 = \G_m$, and that  $S_\infty^G$ is infinite.  If $G$ is not semi-constant,  Proposition 4, combined with Part (A) of the Main Lemma, implies that a multiple of $s$ factors through a translate of $H = \G_m$  by a constant section $s_0 \in \G_m(\C)$, so $s$ factors through a torsion translate of $\G_m$, and $p = \pi(s)$ is torsion, contradicting the hypothesis of Theorem 2. So, $G$ must be semi-constant, $E = E_0 \times S$ must be isoconstant, and the argument just described shows that $s$ must satisfy Conclusion (i) of the Main Lemma, for some $s_0   \in \G_m(\C)$. The mere existence of a Ribet section $s_R := s- s_0$ of $G/S$ implies that $p = \pi(s) = \pi(s_R)$ and $q$ are antisymmetrically related. Moreover, by Theorem 3.i,    $s_R(\overline \lambda)$  is a torsion point on $G_{\overline \lambda}$ whenever $s(\overline \lambda)$ is so, since $\pi(s_R(\overline \lambda)) = \pi(s (\overline \lambda)) = p(\overline \lambda)$ is then a torsion point of $E_0$. There are infinitely many such $\overline \lambda$'s, so at least one. Consequently the  constant section $s_0 \in \G_m(\C)$ is torsion, and (a multiple of) $s$ is a Ribet section of $G/S$. This concludes the argument for Theorem 2.

 \medskip
\noindent
{\bf Part (B) of the Main Lemma}.

 \medskip

Just as for the Main Theorem (cf. Theorem 3), let us right now deal with the ``if" side of Part B of the Main Lemma.  

\smallskip

The periods $\Pi_G$ of $G$ are defined over  the subfield $F_q$ of $F_{pq}$ (see \S 5.1, or the explicit formula given in Appendix I, \S 7.1, or the footnote below), so, clearly, $log_G(s)$ lies in $F_{pq}$ if $s- s_0$ is a torsion section. When $s- s_0 := s_R$ is a Ribet section, an explicit formula for $log_G(s_R)$ in terms of $u$ and $\zeta_\lambda(u)$ is given in Appendix I, \S 7.2, from which the rationality of $log_G(s_R)$ over $F_{pq}$ immediately follows. In fact, we will prove this in a style closer to Manin-Mumford issues in Lemma 3 of \S 6 below\footnote
{~ More intrinsically, concerning the field of definition of $\Pi_G$ : the Cartier dual of the one-motive $[0 \rightarrow G]$ is the one-motive $[\Z \rightarrow \hat E]$ attached to  $q \in \hat E(S)$, so their fields of (generalized) periods coincide, and $K(\Pi_G) = F^{(1)}_G  \subset F^{(2)}_G = F_q$. Concerning $log_G(s_R)$ :  in the notations of \S 1.2,  it suffices to consider   the generic  Ribet section $s_R$  of the semi-abelian scheme ${\cal P}_0$, viewed as an extension ${\cal G}_0$ of $E_0$ by $\G_m$, over the base $\hat E_0$. As mentioned there, its image $W_R$  is a special curve of the mixed Shimura variety ${\cal P}_0$. Therefore, the inverse image of $W_R$ in the uniformizing space of ${\cal P}_0$ is an algebraic curve. In the notations of \S 6,   the statement amounts to the vanishing of $\tau_{s_R}$, and could alternatively be deduced from the self-duality of the one-motive $[M_{s_R} : \Z \rightarrow G]$ attached to the Ribet section, cf. \cite{B-E}.}. The last case forces $G$ to be an isotrivial extension. In the notations of \S 4, we then have $s - s_0 = (\delta, p)
 \in G(S)$,  with $\delta$ a root of unity, so $log_G(s)$ is rational over the field $F_p$.  
 
 \smallskip
 
 As for  the ``only if" side of Part (B) not covered by Part (A), we must show that if (a multiple by a positive integer  of) $s- s_0$ is a non constant section  $\delta$ of $\G_m(S)$, then $\ell := log_{\G_m} (\delta)$ is transcendental over $F_{pq}$. But then, $p -  \pi(s_0)$ is a torsion section of $E/S$, so $F_{pq} = F_q$, and the statement follows  from Lemma 1 of \S 4, with $q$ playing the role of $p$.  
 
 \medskip
 In conclusion, we have reduced the proof of the Main Theorem (more specifically, of Theorems 1 and 2) to defining the field $F_{pq}$,  proving Proposition 4, and   proving Part (A) of  the Main Lemma.

\section{A warm up : the case of direct products }

In this Section, we perform the above-mentioned tasks under the assumption that $G$ is an isotrivial extension,  thereby establishing
Theorem 1, as stated in \S 2. Without loss of generality, we assume that $G = \G_m \times E$, i.e. $q = 0$  (so, the field $F_{pq} = F_{p0}$ will coincide with $F_p$). Of course, if $E/S$ is isoconstant,  say $E = E_0 \times S$, then $G = G_0 \times S$ with $G_0 = \G_m \times E_0/\Q^{alg}$, and Theorem 1 follows from  Hindry's theorem \cite{Hi}; in this isoconstant case, the  strategy we are here following reduces to that of \cite{PZ}. 

\medskip

We first rewrite in concrete terms  the logarithms $U, u$ and their Betti presentations, under no assumption on the elliptic scheme $E/S$ nor on its section $p$. We fix global differential forms\footnote
{~ÊWhen the modular invariant $j(\lambda)$  is constant, i.e. when $E/S$ is isoconstant, we tacitly assume that $E = E_0 \times S$, with $E_0/\C$, and that the chosen differential of first and second kind $\omega, \eta$ are {\it constant} (i.e. come from $E_0/\C$). In particular, the periods $\omega_1, \omega_2$ and quasi-periods $\eta_1, \eta_2$ are constant. The Weierstrass functions are those of $E_0$, and we drop the index  $\lambda$ from their notation.}
 of the first and second kind $\omega, \eta$ for $E/S$, and for any $\lambda \in \Lambda$, we let 
 $$\wp_\lambda~,~ \zeta_\lambda~,~ \sigma_\lambda$$
  be the usual Weierstrass functions attached to the elliptic curve $E_\lambda/\C$ and its differential forms $\omega_\lambda, \eta_\lambda$. We also fix an elliptic logarithm of the point $p({\overline \lambda}_0)$, and extend it   to an analytic function $u(\lambda) = log_E(p(\lambda)) = Arg \wp_\lambda(p(\lambda))$ on $\Lambda$.  Similarly, we fix a basis of periods and quasi-periods for $E_{{\overline \lambda}_0}$, and extend them   to analytic functions $\omega_1(\lambda), \omega_2(\lambda), \eta_1(\lambda), \eta_2(\lambda)$ (of hypergeometric type if $E/S$ is the Legendre curve). There then exist uniquely defined real-analytic functions $b_1, b_2$ with values in $\R$ such that
$$\forall \lambda \in \Lambda, u(\lambda) = b_1(\lambda) \omega_1(\lambda) + b_2(\lambda) \omega_2(\lambda), \qquad \qquad ({\frak R}_u)$$
and the Betti presentation of $log_E(p(\lambda)$ is given by 
$$u^B(\lambda) := log_E^B(p(\lambda)) = (b_1(\lambda), b_2(\lambda)) \in \R^2. $$

\medskip
We now go to $G  = \G_m \times E$ over $\Lambda$. The section $s : \Lambda \rightarrow G $ has two components $(\delta, p)$, where $\delta : S \rightarrow \G_{m/S}$ is expressed by a rational function on $S$. We fix a classical logarithm of $\delta({\overline \lambda}_0)$ and extend it  to an analytic function $\ell (\lambda) := log_{\G_m}(\delta(\lambda))$  on $\Lambda$.  With these notations, the section $log_G \circ s$ of $(Lie G^{an})/\Lambda$ is represented by the analytic map
$$ \Lambda \ni \lambda \mapsto  log_G(s(\lambda)):= U(\lambda) =  \left( \begin{array}{c} 
 \ell(\lambda)\\ u(\lambda)  \end{array}   \right)  \in \C^2 = (Lie G)_\lambda.$$
The $\Z$-local system of periods $\Pi_G$ admits the basis 
$$\varpi_0(\lambda) =   \left( \begin{array}{c} 
  2 \pi i \\ 0  \end{array}   \right), \varpi_1(\lambda) =   \left( \begin{array}{c} 
  0 \\  \omega_1(\lambda)  \end{array}   \right), \varpi_2(\lambda) =   \left( \begin{array}{c} 
  0 \\  \omega_2(\lambda) \end{array}   \right), $$ 
  and the Betti presentation   of  $log_G(s(\lambda))$ is given by
  $$\Lambda \ni \lambda \mapsto U^B(\lambda) := log_G^B(s(\lambda))  = (a(\lambda), b_1(\lambda), b_2(\lambda)) \in \C \times \R^2 = \R^4,$$
  where $a, b_1, b_2$ are the unique real analytic functions on $\Lambda$ satisfying  
    $$\forall \lambda \in \Lambda, \left( \begin{array}{c} 
 \ell(\lambda)\\ u(\lambda)  \end{array}   \right) =   a(\lambda) \left( \begin{array}{c} 
  2 \pi i \\ 0  \end{array}   \right) + b_1(\lambda)  \left( \begin{array}{c}    0\\   \omega_1(\lambda)  \end{array}   \right) + b_2(\lambda)  \left( \begin{array}{c} 0\\   \omega_2(\lambda)    \end{array}   \right). \quad ({\frak R}_{\ell, u})$$
 We then set ${\cal S} = U^B(\Lambda) \subset \R^4$ as usual. 
 
 \smallskip
 Since $q = 0$, the tower of function field extensions of $K = \C(\lambda)$  to be considered  take here the following simple forms :
 \[\begin{array}{ccccccc}
          F^{(1)}   =    K(\omega_1, \omega_2) & ,  & F^{(1)}_{p0} := F^{(1)}_p  =     F^{(1)}(u)  =  K(\omega_1, \omega_2, u)   ,     \end{array} \]
 while their differential  extensions         
          \[\begin{array}{ccccccc}
          F^{(2)} = F^{(1)}(\eta_1, \eta_2)  &  ,  &  F^{(2)}_{p0} := F^{(2)}_p = F^{(1)}_p(\zeta_\lambda(u)) =   K(\omega_1, \omega_2, \eta_1, \eta_2, u, \zeta_\lambda(u))               \end{array} \]
involve the Weierstrass $\zeta$ function, and can be rewritten as 
\[\begin{array}{ccccccc}
          F   :=  F^{(2)} =   K(\omega_1, \omega_2, \eta_1, \eta_2) & ,  & F_p := F^{(2)}_{p}  = F( u, \zeta_\lambda(u)) .\end{array} \]
          We point out that since it contains the field of definition $F$ of the periods of $\omega$ and $\eta$, the field $F_p$ depends only on the section $p$, not on the choice of its logarithm $u$, so, the notation is justified. Furthermore, let $\alpha \in {\cal O} = End(E)$ be a non-zero endomorphism of $E$, and let $p_0 \in {\cal E}_0(\C)$ be a constant section of $E$. Then, the section $p' = \alpha p + p_0$ yields the same field $F_{p'} = F_p$ as $p$. In particular, $F_p = F$ if $p$ is a torsion or a constant section of $E/S$.
           
 \subsubsection*{Proof of Proposition 4 when $q = 0$.}
 
  Suppose for a contradiction that $\cal S$ contains a real semi-algebraic curve $\cal C$, and denote by $\Gamma \subset \Lambda$ the inverse image of $\cal C$ in $\Lambda$ under the map $U^B$ (all we will need is that $\Gamma$ has an accumulation point inside $\Lambda$ , but  it is in fact a real curve). We are going to study the restrictions to   $\Gamma$ of the functions
 $$a, b_1, b_2, u, \ell, \omega_1, \omega_2.$$
 Recall that all these are functions of $\lambda \in \Lambda$. In view of the defining relation $({\frak R}_{\ell, u})$, the transcendence degree of the functions $u, \ell$ over the field $\C(\omega_1, \omega_2, a, b_1, b_2)$ is at most 0. When restricted to $\Gamma$, the latter field has transcendence degree $\leq 1$ over $\C(\omega_1, \omega_2)$, since $U^B(\Gamma) = (a, b_1, b_2)(\Gamma)$ is the algebraic curve $\cal C$. So, the restrictions to $\Gamma$ of the two functions $u, \ell$ generate over $\C( \omega_{1|\Gamma}, \omega_{2|\Gamma})$ a field of transcendence degree $\leq 1 + 0= 1$, and are therefore algebraically dependent over $\C( \omega_{1|\Gamma}, \omega_{2|\Gamma})$. Since $\Gamma$ is a real curve of the complex domain $\Lambda$, the complex analytic functions $u, \ell$ are still algebraically dependent over the field of $\Lambda$-meromorphic functions $\C(\omega_1, \omega_2)$, i.e.
 $$tr. deg._{\C(\omega_1(\lambda), \omega_2(\lambda))} \;  \C\big(\omega_1(\lambda), \omega_2(\lambda),  u(\lambda), \ell(\lambda)\big) \leq 1. $$
  Now, assume as in Proposition 4 that $p$ has infinite order, and if $E/S$ is isotrivial, that $p$ is not constant. Then, Andr\'e's Theorem 3 in \cite{An} (see also \cite{BT}, Thm. 5) implies that $u(\lambda)$ is transcendental over the field $F^{(1)} = K \big(\omega_1(\lambda), \omega_2(\lambda))$. The previous inequality therefore says that the function $\ell(\lambda)$ is algebraic over the field $F^{(1)}_{p}= K\big(\omega_1(\lambda), \omega_2(\lambda),  u(\lambda)\big)$, or equivalently, that the field of definition of   $log_G(s(\lambda)) = (\ell(\lambda), u(\lambda))$ is algebraic over $F^{(1)}_{p0}$, hence also over $F_{p0} = F^{(2)}_{p0}$, and Proposition 4 is proved when $q = 0$.

\subsubsection*{Proof of Main Lemma (A) when $q  = 0$}
Before giving this proof, let us point out that the advantage of the fields $F = F^{(2)}, F_{pq} =  F^{(2)}_{pq}$ over their first kind analogues is that they are  closed under the derivation $' = \partial /\partial \lambda$.  Moreover, by Picard-Fuchs theory, they are Picard-Vessiot (i.e. differential Galois) extensions of $K$. Since $\ell(\lambda)$ satisfies a $K$-rational DE of order 1, $K(\ell)$ and $F_p(\ell) = F_p(log_G(s))$ too are Picard-Vessiot extension of $K$.

 \begin{Lemma} Let $\Lambda$ be a ball in $\C$, let $\{\wp_\lambda, \lambda \in \Lambda\}$ be  a family of Weierstrass functions, with invariants $g_2, g_3$ algebraic over $\C(\lambda)$ and periods $\omega_1, \omega_2$, and let $u$ be an analytic function on $\Lambda$ such that $u, \omega_1, \omega_2$ are linearly independent over $\Q$. If $j(\lambda)$ is constant, we assume that  $g_2, g_3$ too are constant, and that  $u$ is not constant. Let further $\ell$ be a {\rm non-constant} analytic function on $\Lambda$. Assume that $\wp_\lambda(u(\lambda))$ and $e^{\ell(\lambda)} := \delta(\lambda)$ are algebraic functions of $\lambda$, and consider the tower of differential fields $K \subset F \subset F_p$, where $K= \C(\lambda),  F = K(\omega_1, \omega_2, \eta_1, \eta_2), F_p = F(u, \zeta_\lambda(u))$.   Then,
 $$tr.deg._{F} F_p(\ell(\lambda) \big) =  3.$$
 In particular,  $\ell(\lambda)$ is transcendental over $F_p$, i.e. Part (A) of the Main Lemma holds true when $q = 0$.
 \end{Lemma}
This last statement  is indeed equivalent to Part (A) of the Main Lemma when $G \simeq \G_m \times E$ is a trivial extension (or more generally, an isotrivial one), with constant part $G_0 = \G_m$, resp. $G_0 \times S = G$, if $E$ is not, resp. is, isoconstant.  Indeed, with $s = exp_G(\ell, u) = (\delta, p)$ as above, we then have $F_p = F_{pq}$, and $F_p(\ell) = F_p(log_G(s))$. Lemma 1 then says  that  if $log_G(s)$ is algebraic over $F_{pq}$, then either $p$ is a torsion point (so, a multiple of $s$ factors through $\G_m$), or  $E$ isoconstant and  $p = p_0$ is constant (so, the constant section $s_0 = (1, p_0) \in G_0(\C)$ satisfies $s- s_0 \in \G_m(S)$), or $\delta = \delta_0$ is constant (so, the constant section $s_0 = (\delta_0, 0) \in G_0(\C)$ satisfies $s- s_0 \in E(S)$). In all cases, we therefore derive Conclusion (ii) of the Main Lemma.

 \medskip
 \noindent
 {\it Proof}  (of  Lemma 1): this is essentially due to Andr\'e, cf. \cite{An}, but not fully stated there (nor in \cite{BT}). It is proven in full generality in \cite{BM}, but one must look at the formula on top of p. 2786 to see that $K$ can be replaced by $F$ in  Theorem L ... So, it is worth giving a direct proof.
 
 \smallskip
 
  We first treat the case when $E/S$ is not isoconsant. By Picard-Lefchetz, the Picard-Vessiot extension $F = K(\omega_1, \omega_2, \eta_1, \eta_2)$   of $K = \C(\lambda)$ has Galois group $SL_2$. By \cite{BT}, the Galois group of $F_p   = F(u(\lambda), \zeta_\lambda(u(\lambda))$ over $F$ is a vector group ${\cal V}Ê\simeq \C^2$ of dimension 2 (i.e. these two functions are algebraically independent over $F$), while the Galois group of $K( \ell(\lambda))$ over $K = \C(\lambda)$ is $\C$. Since $\C$ is not a quotient of $SL_2$, the Galois group of $F(\ell(\lambda))$ over $F$ is again $\C$. Now, $SL_2$ acts on the former ${\cal V} = \C^2$ by its standard representation, and on the latter $\C$ via the trivial representation, so  the Galois group of $F(u(\lambda), \zeta_\lambda(u(\lambda), \ell(\lambda))$ over $F$ is a subrepresentation $\cal W$ of $SL_2$ in $\C^2 \oplus \C$ projecting onto both factors. Since the standard and the trivial representations are irreducible and non isomorphic, we must have ${\cal W} = \C^2 \oplus \C$. Therefore, 
 $$tr.deg._{F} \; F(u(\lambda), \zeta_\lambda(u(\lambda)),  \ell(\lambda)) =  dim {\cal W} = 3. $$
 
 \smallskip
 We now turn to the case of an isoconstant $E = E_0 \times S$. The field of periods $F$ then reduces to $K$, and $SL_2$ disappears. But since the ambient group $G = \G_m \times  E$ is now isoconstant, we can appeal to Ax's theorem on the functional version of the Schanuel conjecture.  More precisely, since the result we  stated involves the $\zeta$-function, we appeal to its complement on vectorial extensions, see \cite{BK}, Thm. 2.(iii),  or more generally, \cite{Be}, Proposition 1.b,  which implies that  $tr.deg_K K(u, \zeta(u), \ell) = 3$ as soon as $u$ and $\ell$ are not constant. (Actually, Andr\'e's method still applies to the isoconstant case, but requires a deeper argument, involving  Mumford-Tate groups : cf. \cite{An}, Theorem 1, and \cite{BM}, \S 8.2). 
 
 \medskip
 \noindent
 {\bf Remark 2} :  ~(i) Concerning  the proof of Theorem 1 :  in fact, as pointed out by the fourth author in \cite{Za}, p. 79, Comment (v),   since we are here dealing with a direct product, the torsion points yield torsion points on $\G_m$, which lie on the unit circle, a real-curve. So,  in the argument of \S 3.2, the dimension decreases by 1 a priori, and  Bombieri-Pila (real-curves) rather than Pila (real-surfaces) as in Proposition 2 would suffice.

 \smallskip
 
 ~  (ii) As shown by the proof of Proposition 4 ($q$ = 0) above,  it would have sufficed to prove that the (non constant) logarithm $\ell$ is transcendental over the field $\C(\omega_1, \omega_2, u)$. Adjoining $\lambda$ leads to  $F^{(1)}_p$ as base field, and as already said, differential algebra then forces to consider $F^{(2)}_p$.  For a broader perspective on these extensions of the base field, see \S 5.3 below.

\smallskip

~(iii) For the last statement of Lemma 1 to hold, the only necessary  hypothesis is that $\ell$ be non-constant. Indeed, if $p$ is torsion or constant, then $F_p = F$, and $u$ plays no role. But we have preferred to present Lemma 1 and its proof in this way,   as an introduction to the general proofs of \S 5 and \S 6.

\smallskip
 
~ (iv) Conversely, let $q$ be any (non-necessarily torsion or constant)  section of $E/S$, and set $v = log_E(q)$, $F_q = F^{(2)}(v, \zeta_\lambda(v))$ and $F_{pq} = F_p.F_q$, as will be done in \S 5.  The same proof as above shows that
$$\forall  p, q \in E(S), \forall \delta \in \G_m(S), \delta \notin \G_m(\C)), ~\ell := log_{\G_m}(\delta) ~
is  ~ transcendental ~ over ~  F_{pq}.$$
 Indeed, the only new case is when $p$ and $q$ are linealy independent over $End(E)$ modulo the constant part of $E/S$. From the same references and argument as above, replacing ${\cal V} = \C^2$ by ${\cal V}Ê\oplus {\cal V}$ , we deduce that    the transcendence degree of $F_{pq}(\ell)$ over $F$ is equal to 5\;, yielding the desired conclusion on the transcendency of $\ell$.

\subsubsection*{A characterization of Ribet sections}
We close this section on isotrivial extensions by a corollary to Theorem 1,  which plays a useful role in checking the compatibility of the various definitions of Ribet sections : see the equalities $\beta_R =  \beta_J$ in \cite{BE}, and $s_{\tilde R} = s_R$ in Appendix I, Proposition 5.(i) below.

\begin{Corollary}   Let $G/S$ be an extension of $E/S$ by $\G_m$, and let $p$ be a section of $E/S$ of infinite order and not constant if $E/S$ is isoconstant (equivalently, by  Proposition 1.iv, such   that the set $S_\infty^E$ attached to $p$ is infinite). Let further $s^\dagger$ and $s$ be two sections of $G/S$ such that $\pi \circ s^\dagger  = \pi \circ s = p$. Assume that for all but finitely many (resp. for infinitely  many) values of $\overline \lambda$ in $S_\infty^E$, the point $s^\dagger (\overline \lambda)$ (resp. $s(\overline \lambda)$) lies in $G_{tor}$, i.e. that ~Ê{\rm   $s^\dagger$ (resp. $s$) ``lifts almost all (resp. infinitely many) torsion values of $p$ to torsion points of $G$".}   Then, there exists a torsion section $\delta_0 \in \G_m(\C)$ (i.e. a  root of unity) such that $s = s^\dagger + \delta_0$. 
\end{Corollary}

\noindent
{\it Proof} : Let  $\delta_0 := s- s^\dagger \in \G_m(S)$. We know that $s^\dagger$ lifts almost {\it all} torsion points $p(\overline \lambda) \in E_{tor}$ to points in $G_{tor}$. If  $s$ does so for infinitely many of them, then, so does the section $s_1 := (\delta_0, p)$ of the direct  product $G_1 = \G_m \times E/S$. This contradicts Theorem 1, unless the projection  $\delta_0(S)$ of $s_1(S)$ to $\G_m$ is reduced to a root of unity. 

\bigskip

 It is interesting to note that in this way, Theorem 1 on the trivial extension $G_1$ has an impact on extensions $G$  which need not at all be isotrivial. For instance, if $G$ is semi-constant and $E_0$ has CM, Corollary 3, applied to the Ribet section $s^\dagger = s_R$, shows that up to isogenies, $s_R$ is the only section which lifts infinitely many torsion values of $\pi(s_R)$
  to torsion points of $G$. We also point out that since   the elliptic scheme $E \simeq E_0 \times S$ is here constant,  Hindry's theorem  on the constant semi-abelian variety $\G_m \times E_0$ suffices to derive this conclusion.
 
 \medskip

On the other hand, assume that $G$ is an isotrivial extension. Then, there exists a subgroup scheme $E^\dagger/S$ of $G/S$ such that the restriction  $\pi^\dagger$ of $\pi : G \rightarrow E$ to $E^\dagger$ is an $S$-isogeny.   Any section $s^\dagger$ of $G/S$ a non zero multiple of which factors through $E^\dagger$ then satisfies the lifting property of the corollary, since $p(\overline \lambda) := \pi^\dagger  \circ s^\dagger (\overline \lambda)$ is a torsion point of $E_{\overline \lambda}$ if and only if  $s^\dagger (\overline \lambda)$ is a torsion point of  $G_{\overline \lambda}$ . By Corollary 3, such sections $s^\dagger$ are, up to a root of unity, the only section $s$ above $p = \pi \circ s^\dagger $ for which  $s(S) \cap G_{tor}$ is infinite. Of course, this is  (after an isogeny) just a rephrasing of Theorem  1, but it shows the analogy between these  ``obvious" sections and the Ribet sections. This is a reflection of the  list of special curves of the mixed Shimura variety described in \S 1.2.

\section{The general case}

\subsection{Fields of periods and the Main Lemma}

Apart from the statement of Lemma 2 below, we henceforth make no assumption on the extension $G$ of $E/S$ by $\G_m$. So the section $q \in \hat E(S)$ which parametrizes $G$ is arbitrary. Concerning the elliptic scheme $E/S$, we  recall the notations of \S 4, and in particular, the fields of periods $F^{(1)} = K(\omega_1, \omega_2)$  of $E$ and its differential extension   $F := F^{(2)} = F^{(1)}(\eta_1, \eta_2)$.  We identify $\hat E$ and $E$ in the usual fashion, 	and  denote by $v = log_E(q)$  a logarithm of the section $q$ over $\Lambda$. We recall that the field $F_q := F^{(2)}_q = F^{(2)}(v, \zeta_\lambda(v))$   depends only on $q$, and coincides with $F^{(2)} = F$ when $q$ is a torsion section, i.e. when $G$ is isotrivial.   We will also use the notations of \S 3.3 on the local system of periods  $\Pi_G = \Z \varpi_0 \oplus \Z \varpi_1 \oplus  \Z \varpi_2 \subset Lie G^{an}(\Lambda)$ of $G^{an}/\Lambda$.

\medskip

Consider the extension $F^{(1)}_G = K(\varpi_0,  \varpi_1, \varpi_2) = K(\varpi_1, \varpi_2)$  of $K = \C(\lambda)$ generated by  the elements  of $\Pi_G$. Since $\Pi_G$  projects onto $\Pi_E$ under $d\pi$, whose kernel has relative dimension 1,  this field is an extension of $F^{(1)} = K(\omega_1, \omega_2)$ of transcendence degree $\leq 2$. So, the field $F^{(2)}_G$ generated by $\Pi_G$ over $F := F^{(2)}$ has transcendence degree $\leq 2$. In fact, the duality argument mentioned in Footnote (3), or more explicitly, the computation given in Appendix I, shows that  
$$F^{(2)}_G =   F^{(2)}(v, \zeta_\lambda(v)) := F_q ,$$
i.e $F^{(2)}_G$ coincides with the differential extension $F_q$ of $F^{(2)} = F $ attached to $q$. So,  $F^{(2)}_G$ is in fact a Picard-Vessiot extension of $K$.

\medskip
A more intrinsic way to describe these ``fields of the second kind" is to introduce the {\bf universal vectorial extension} $\tilde E/S$ of $E/S$, cf. \cite{BP}. This   is an $S$-extension of $E/S$ by the additive group $\G_a$, whose local system of periods $\Pi_{\tilde E}$ generates the field $F^{(2)}$.  The universal vectorial extension $\tilde G/S$ of $G/S$ is the fiber product $G \times_E \tilde E$, and its local system of periods $\Pi_{\tilde G}$ generates the field $F^{(2)}_G$. Now, for both $\tilde E$ and $\tilde G$ (and contrary to $E$ and $G$), these local systems generate the spaces of horizontal vectors of   connections $\partial_{Lie \tilde E}, \partial_{Lie \tilde G}$ on $Lie \tilde E/S, Lie \tilde G/S$. This explains why the fields  $K(\Pi_{\tilde E}) = F^{(2)} = F$ and $K(\Pi_{\tilde G}) = F^{(2)}_G = F_q$ are  Picard-Vessiot extensions of $K$.

\medskip
Let now $s$ be a section of $G/S$, and let $U = log_G(s) \in Lie G^{an}(\Lambda)$ be a logarithm of $s$ over $\Lambda$. As usual, set $p = \pi(s) \in E(S), u = log_E(p) = d\pi(U) \in Lie E^{an}(\Lambda)$. Since  $Ker (d\pi)$ has relative dimension 1,   the field generated over $K$ by $log_G(s)$ is an extension of $K(log_E(p))$ of transcendence degree $\leq 1$, so $log_G(s)$ has transcendence degree $\leq 1$ over  $F_p = F(u, \zeta_\lambda(u))$. Finally, set
$$F_{pq} := F_p.F_q  ~, ~L = L_s :=  F_{pq}(log_G(s)).$$
The  field $F_{pq} = F^{(2)}_{pq}$ is the field of generalized periods of $\{E, p,q\}$ promised in \S 3.3.  Since it contains $F^{(2)}_G$,  the   field $L = L_s$  depends only on $s$, not on the choice of its logarithm $U = log_G(s)$, and  is an  extension of $F_{pq}$, of transcendence degree $\leq 1$.  In fact, the explicit formulae of Appendix I show that for $q \neq 0$ and $p \neq 0, - q$ :
$$L_s = F_{pq}\big(\ell_s - g_\lambda(u, v)  \big), ~Ê{\rm where}~ g_\lambda(u, v) =    log_{\G_m} \frac{\sigma_\lambda(v+u)}{\sigma_\lambda(v)  \sigma_\lambda(u)} $$
is a Green function attached to the sections $\{p, q\}$ of $E/S$, and $\ell_s = log_{\G_m}(\delta_s)$ for some rational function $\delta_s \in K^*$ attached to the section $s$ of $G/S$.
This formula implies that $L$ is a differential field, but $L$ is even a Picard-Vessiot extension of $K$. One way to check this is to relate $log_G(s)$ to an integral of a differential of the 3rd kind on $E$, with integer, hence {\it constant} residues, and to differentiate under the integral sign. Another way consists in lifting $s$ to a section $\tilde s$ of $\tilde G/S$, projecting to $\tilde p \in \tilde E(S)$. Then, for any choices $\tilde u = log_{\tilde E} (\tilde p)$ and $\tilde U = log_{\tilde G}(\tilde s)$ of   logarithms of $\tilde p$ and $\tilde s$, the field $F^{(2)}(\tilde u) = F(u, \zeta_\lambda(u)) = F_p$ is contained in $F_q(\tilde U)$, which can therefore be written as $F_{pq}(\tilde U)$, and the latter field  $F_{pq}(\tilde U)$ coincides with $F_{pq}(U) = L$, since $\tilde U$ lifts $U$ and $\tilde u$ in the fiber product $Lie \tilde G = Lie G \times_{Lie E} Lie \tilde E$. Now, in the notations of \cite {BP}, \cite {BM}, $\tilde U$ is a solution of the inhomogenous linear system $\partial_{Lie \tilde G}(Ê\tilde U) = \partial \ell n_{\tilde G}( \tilde s)$, which, on the one hand, is defined over $K$, and on the other hand, admits  $F_q(\tilde U)$  as field of solutions. So, $L = F_{q}(\tilde U)$ is indeed a Picard-Vessiot extension of $K$. By the same argument, applied to the differential equation  $\partial_{Lie \tilde E}(Ê\tilde u) = \partial \ell n_{\tilde E}( \tilde p)$, we see anew that $F_p$ is a Picard-Vessiot extension of $K$. Notice, on the other hand, that $F_q(U)$ is in general not a differential extension of $F_q$ (it contains $u$, but not $\zeta_\lambda(u)$). 

\medskip
The following diagram  summarizes these notations (and proposes other natural ones...)~:
 $$ \begin{array}{ccccc} 
&& L & \\  && \uparrow &  & \\ && F_{pq} && \\ &\nearrow & & \nwarrow &\\ F_q &&&& ~F_p  \\ &   \nwarrow & &   \nearrow & \\&& F && \\   &&  \uparrow  & \   & \\ && K && \end{array}
\qquad \begin{array}{ccccc} 
  L = F_{pq}(log_G(s)) = F_{\tilde G}(log_{\tilde G}(\tilde s))  = F_{pq}\big( \ell_s - g(u,v) \big)  \\  &&  &  & \\   F_{pq}= F_p.F_q = F(u, \zeta(u), v, \zeta(v))   \\ &  & &  &\\   F_q = F^{(2)}_G := F_{\tilde G} = F_{\tilde E}(log_{\tilde E}(\tilde q))  \qquad  F_p  =F_{\tilde E}(log_{\tilde E}(\tilde p))  \\ &    & &    & \\ F = F^{(2)} := F_{\tilde E} = K(\omega_1, \omega_2, \eta_1, \eta_2)\\   &&     &   & \\   K = \C(\lambda)  \end{array}$$

All the notations of the Main Lemma  have now been discussed,  and  we can restate its Part (A) in the non-isotrivial case, as follows : 

\begin{Lemma} {\rm ($=$  Main Lemma for $q$ non-torsion)}Ê~  {\it With $S/\C$, let $G/S$ be a {\rm non isotrivial}  extension by $\G_m$ of an elliptic scheme $E/S$, parametrized by a section $q$ of $\hat E/S$, and  let $G_0$ be the constant part of $G$. Let further $s $ be a section of $G/S$, with projection $p = \pi \circ s $ to $E/S$,  and let $F_{pq} = F_p.F_q  \supset F$ be the   field  of generalized periods of $\{E, p, q\}$.   Assume that  $log_G(s)$ is algebraic over $F_{pq}$. Then, there exists a constant section $s_0 \in G_0(\C)$ such that

\smallskip

i) either  $s - s_0$ is a Ribet section;

ii) or  $s - s_0$ is a torsion section.

\smallskip
\noindent
In other words, if $s$ is not a constant translate of a Ribet or a torsion section of $G/S$, then $ g_\lambda(u,v) - \ell_s$ is transcendental over $F_{pq}$.
}
\end{Lemma}
Conclusion (ii) of Lemma 2 appears to be stronger than Conclusion (ii) of  Part (A) of the Main Lemma, but is in fact equivalent to it when $G$ is not isotrivial. Indeed, in this case,  $\G_{m/S}$ is the only connected strict subgroup scheme of $G$. Now, if a multiple by a non-zero integer $N$ of the section $s- s_0$ factors through $\G_{m}$, i.e. is of the form $\delta$ for some section $\delta \in \G_m(S)$,  then, $p = \pi(s)$ is a torsion or a constant section of $E/S$, so $F_{pq} = F_q$, and $F_{pq}(\log_G(s)) = F_q(\ell)$, where $\ell = log_{\G_m}(\delta)$. By assumption, $\ell$ is then algebraic over $F_q$, and Lemma 1 of \S 4  implies that $\delta = \delta_0 \in \G_m(\C)$ is a constant section. Considering the constant section $s'_0 = s_0 - \frac{1}{N}  \delta_0$ of $G/S$, we derive   that $s -s '_0$ is a torsion section, i.e. that Conclusion (ii) of Lemma 2 is fulfilled. 

\subsection{Reducing the Main Theorem to the Main Lemma}

 In view of \S 4, we could now restrict to the case of a non isotrivial extension $G$, prove Proposition 4 in this case, and finally prove Lemma 2, thereby concluding the proof of Theorem 2. However, as announced above, we will  remain in the general case, and make no assumption on $q$. 
   
 \smallskip
 We now prove Proposition 4,  extending the  pattern of proof of \S 4 to the general case. We recall the notations of Proposition 4, including the fundamental assumption that $p = \pi(s)$ is neither a torsion nor a constant section of $E/S$.  By Lemma 1 , this condition implies that $F_p$ has transcendence degree $2$ over $F$, hence that {\it $u(\lambda) = log_E(p(\lambda))$ is transcendental over the field $F$.}

 \smallskip
 Consider the following tower of fields of functions on $\Lambda$, where the lower left (resp. upper right) ones are generated by complex (resp. real) analytic functions. The inclusions which the NE-arrows represent come from the definition of $a, b_1, b_2$ in terms of $log_E(p(\lambda)) = u(\lambda), log_G(s(\lambda)) = U(\lambda)$, cf. Relations  $({\frak R}_u), ({\frak R}_U)$ of \S 3.3; the inclusions of the NW-arrows on the left come from the definition of the fields  of periods $F, F_{pq}$; those of the  NW-arrows on the right are obvious.  
 
 $$ \begin{array}{rcccccccc} 
  & &    F_{pq}(a, b_1, b_2) &&&& ({\frak R}_U):&log_Gs = a \varpi_0 + b_1 \varpi_1 + b_2 \varpi_2    \\  & \nearrow & &   \nwarrow && &   \\ &&&&  &  &&  \\  F_{pq}(log_G(s)) &&&& F_{pq}(b_1, b_2)  & & \\ &&&&&&&  \\  & \nwarrow   && \nearrow &&     \nwarrow &   \\ & &F_{pq}   &&&&  F(b_1, b_2)  \\ &&&&&&&  \\ & & &   \nwarrow &&   \nearrow  && \\ &&&&  F(log_E(p)) &&({\frak R}_u):& u :=  log_E p = b_1 \omega_1 + b_2 \omega_2 \\   &&&&   \uparrow  &&  \\ &&&&  F && \end{array}. $$

\medskip
   Now, assume for a contradiction that the real surface  $\cal S$ contains a semi-algebraic curve $\cal C$. As in \S 4, consider the real curve   $\Gamma = (U^B)^{-1}({\cal C}) \subset \Lambda  \subset S(\C)$, and denote by a lower index $\Gamma$ the restrictions to $\Gamma$ of the
 various functions of $\lambda$ appearing above;  similar notation $F_{|\Gamma}, F_{pq|\Gamma}$, etc, for the fields they generate. For instance, since $(a_{|\Gamma}, b_{1|\Gamma}, b_{2|\Gamma})$ parametrize the algebraic curve ${\cal C}$, these three functions generate a field of transcendence degree 1 over $\R$, and $tr. deg._{F_{|\Gamma}} {F_{|\Gamma}}(b_{1|\Gamma}, b_{2|\Gamma} ) \leq 1$. But by the result recalled above and the principle of isolated zeroes, $u_{|\Gamma}$ is transcendental over $F_{|\Gamma}$. Therefore, the restriction to $\Gamma$ of the field $F(b_1, b_2)$ is an algebraic extension of the restriction to $\Gamma$ of the field $F(u)$. We may abbreviate this property by saying that $F(u)$ and $F(b_1 , b_2)$ are essentially equal over $\Gamma$. Going up north-west in  the tower, we deduce that the fields $F_{pq}$ and $F_{pq}(b_1, b_2)$ are essentially equal over $\Gamma$. 
 
 Notice that $b_{1|\Gamma}$ and $ b_{2|\Gamma}$ are not both constant   since ${F_{|\Gamma}}(b_{1|\Gamma}, b_{2|\Gamma} ) := \big(F(b_1, b_2)\big)_{|\Gamma}Ê\supset \big(F(u) \big)_{|\Gamma}$ is a transcendental extension of $F_{|\Gamma}$. So, $a_{|\Gamma}$ must be algebraic over  $\R(b_{1|\Gamma}, b_{2|\Gamma} )$. Therefore, $\big(F_{pq}(a, b_1, b_2)\big)_{|\Gamma}$ is an algebraic extension of $\big(F_{pq}(b_1, b_2) \big)_{|\Gamma}$, hence of the essentially equal field $\big(F_{pq}\big)_{|\Gamma}$ , and we deduce that the intermediate field $\big(F_{pq}(log_G(s)) \big)_{|\Gamma}$ is algebraic over $\big(F_{pq}\big)_{|\Gamma}$. But $log_G(s(\lambda)) = U(\lambda)$ is a complex analytic map, so, by isolated zeroes,  $F_{pq}(log_G(s(\lambda))$  must also be algebraic over $F_{pq}$. This concludes the proof of Proposition 4.

 \subsection{The role of $K$-largeness}
 
 The change of base field from $K$ to $F_{pq}$ can be viewed as the ``logarithmic" equivalent of the passage from $K$ to the field $K_{\bf G}^\sharp$ generated by the Manin kernel ${\bf G}^\sharp := Ker(\partial \ell n_{\bf G})$ of an algebraic $D$-group ${\bf G}/K$, which one encounters in the study of the exponentials of algebraic sections of $Lie({\bf G})/S$, as in \cite{BP}. A Manin kernel has in fact  already appeared in the present paper at the level of the $D$-group $\tilde E$ : in this case, the group of $K^{alg}$-points of $\tilde E^\sharp$ projects onto $E_{tor}Ê\oplus {\cal E}_0(\C)$, and the recurrent hypothesis made on the section $p$ (that it be neither torsion nor constant) exactly means that none of  its   lifts $\tilde p$ to $\tilde E(S)$   lies in $\tilde E^\sharp$. 
 
 \medskip
 For elliptic curves, the field of definition $K_{\tilde E}^\sharp$ of  $\tilde E^\sharp$ is always algebraic over $K$, and one says that the $D$-group $\tilde E$ is $K$-large. This is the hypothesis required on ${\bf G}$  for the Galois theoretic approach to the proof of the relative Lindemann-Weierstrass theorem of \cite{BP}, \S 6  (see also \cite{BM}, \S 8.1  and \S 4.3 of S. \& C. for the abelian case). But    the third author \cite{Pl} has checked that it can be extended to non $K$-large groups (those with $K_{\bf G}^\sharp$ transcendental over $K$), as follows  :

\bigskip
\noindent
{\bf Theorem} : {\it Let ${\bf G}$ be an almost semi-abelian algebraic $D$-group over $K$, let $a \in Lie{\bf G}(K)$, and let $y = exp_{\bf G}(\int a) \in {\bf G}(K^{diff})$ be a solution of the equation $\partial \ell n_{{\bf G}}(y) = a$. Then, $tr.deg.(K^\sharp_{\bf G}(y)/K^\sharp_{\bf G})$ is the smallest among dimensions of connected algebraic $D$-groups $H$ of ${\bf G}$ defined over $K$ such that $a \in Lie H + \partial \ell n_{\bf G}({\bf G}(K))$, equivalently such that $y \in H + {\bf G}(K) + {\bf G}^\sharp$. Moreover, $H^\sharp(K^{diff})$ is the Galois group of $K^\sharp_{\bf G}(y)$ over $K^\sharp_{\bf G}$.}

\bigskip
By $\int a$, we here mean any $x \in Lie{\bf G}(K^{diff})$ such that $\partial_{Lie {\bf G}}x = a$. When $x$ lies in $Lie {\bf G}(K)$, this leads to results of ``Ax-Lindemann" type (as used in \cite{PZ}, \cite{Pi}), whereas the transcendence results required by  the present strategy ($=$ that of \cite{MZ}) concern the equation  $\partial_{Lie {\bf G}}x  = b$, with $b = \partial \ell n_{{\bf G}}(y)$ for some $y \in {\bf G}(K)$. Notice that when  ${\bf G}$ is the universal vectorial extension $ \tilde G$ of our semi-abelian scheme $G = G_q$, the analogous field  $K_{Lie {\bf G}}^\sharp$ is precisely the field of periods $F_q = F_{\tilde G}$ of $\tilde G$; see  \cite{BM}, \S 2 of S \& C, for a justification of this analogy, and Remark 3.ii below for the adjunction of $F_p$ in the base field.  

\medskip
 More generally, given  $(x, y)    \in (Lie {\bf G}  \times {\bf G})$, analytic over a ball $\Lambda \subset S(\C)$ and linked by the relation  $\partial \ell n_{{\bf G}}(y) = \partial_{Lie {\bf G}}x$ (i.e. essentially $y = exp_{\bf G}(x)$), one may wonder which extension of the  Ax-Schanuel theorem holds for the transcendence degree {\it over} $K_{Lie {\bf G}}^\sharp. K_{\bf G}^{\sharp}$ of the point $(x, y)$. The case when $x$ is algebraic over $K_{Lie{\bf G}}^\sharp$ includes the study of Picard-Painlev\'e sections; in this direction, see \cite{Ha}, Lemma 3.4. The case when $x$ has transcendence degree $\leq 1$ over $K_{Lie{\bf G}}^\sharp$  would be of particular interest, as it  occurs when the Betti coordinates of $x$ parametrize an algebraic curve, and this may  pull  the present strategy back into the ``exponential" frame-work of Ax-Lindemann-Weierstrass.  

\medskip
\noindent
{\bf Remark 3}.-  i) Just as in \S 4, all we need to know for relative Manin-Mumford statements is the transcendence degree of $log_G \circ s$ over the   field of generalized periods $F_{pq}$ of $\{E, p, q\}$. More clearly put,  the transcendence degree (i.e. the differential Galois group) of $F_{pq}$ over $K$ plays no role. Of course, $Gal(F_{pq}/K)$ will  come as a help during the proof of the Main Lemma, in parallel with the role of $SL_2$ in \S 4. But in the notations of the diagram of \S 6, we must merely  compute $Gal(L/F_{p,q}) = Im(\tau_s) \subset \C$, and show that under the hypotheses of the Main Lemma,  $\tau_s$ vanishes    only if one of its conclusions (i), (ii) is satisfied.

 \smallskip
 
ii) Adjoining the field $F_p$ to the base field $F_q$ comes in naturally, since the Picard-Vessiot extension $F_q(\tilde U)$ of $K$ automatically contains it. But another advantage of the compositum $F_{pq}$ is that the roles of $p$ and $q$ become symmetric in the statement of the Main Lemma. The explicit formula for $L$ given by the Green function makes this apparent. More intrisically, the field $L_s = F_{pq}(log_G(s))$ is the field of periods of the  smooth one-motive $M = [M_s: \Z \rightarrow G]$ over $S$,  in the sense of \cite{De},  attached to the section $s \in G(S)$, with $p = \pi(s)$ and $q \in \hat E(S)$ parametrizing the extension $G$. By biduality, $p$ parametrizes an extension $G'$ of $\hat E$ by $\G_m$, and $s$ may be viewed as a section $s'$ of $G'$ above the section $q$ of $\hat E(S)$. The Cartier dual of $M$ is the one-motive $M' = [M'_{s'} : \Z \rightarrow G']$ attached to this section $s'$, and its field of periods $F_{M'} = L'_{s'} = F_{qp}(log_{G'}(s'))$ coincides with $F_M$, since they are the Picard-Vessiot extensions of two adjoint differential systems. So, although $log_G$ and $log_{G'}$ have no direct relations, the fields which  $log_G(s)$ and $log_{G'}(s')$ generate {\bf over $F_{pq}$} are the same. Similarly,   the structures of $G$ and $G'$ usually differ a lot,  but the conclusions (i) and (ii) of the Main Lemma turn out to be invariant under this duality. See Case ($\bf SC 2$) of \S 6 for a concrete implementation of this remark.

\smallskip
iii) The above symmetry is best expressed in terms of the Poincar\'e bi-extension $\cal P$, resp. ${\cal P}'$, of $E \times_S \hat E$, resp. $\hat E \times_S E$, by $\G_m$. As recalled in the Introduction (cf. \cite{De}), a section $s$ of $G = G_q$ above $p$ corresponds to a trivialization  of the $\G_m$-torsor $(p, q)^*{\cal P} \simeq (q, p)^*{\cal P}'$. Then, the inverse image $\varsigma$ of this trivialization under the uniformizing map 
$$ \C^3 \times \tilde S  \rightarrow {\cal P}^{an}$$
 of $\cal P$ generates $L$ over $F_{pq}$ (here, $\tilde S$ denotes the universal cover of $S^{an}$, for instance the Poincar\'e half-plane when $S = X$ is a modular curve). This view-point turns the  Main Lemma   into a statement  about the transcendency of $\varsigma$ over $F_{pq}$, and explain why  the various types of special curves of the mixed Shimura variety ${\cal P}/X$, as encountered in \S 1.2, occur in its conclusion. 
 
 \smallskip
 iv) ({\it Autocritique} on differential extensions) It would be interesting to pursue the study of this uniformizing map further, as it may lead to a simplification of the present proof of the Main Theorem, where the appeal to differential extensions would be replaced by an Ax-Lindemann statement, extending the recent results of Ullmo-Yafaev \cite{UY} and Pila-Tsimerman \cite{PT} to mixed Shimura varieties. See \cite{Sc} for a perspective on both  approaches.

\section{Proof of the Main Lemma}

 The arguments will    be of the same nature as   in  \S 4, appealing to Ax type results for constant groups, and to representation theory otherwise. As mentioned before the enunciation of the Main Lemma in \S 3.3, similar results appear in  \cite{Be}, Propositions 4.a, 4.b and Theorem 2, but it seems better to gather them here into a full proof.

\medskip
Consider the tower of Picard-Vessiot extensions drawn on the  left part of the following picture : 
 $$ \begin{array}{ccccc} 
&& L && \\  && \uparrow &  & \\ && F_{pq} && \\ &\nearrow & & \nwarrow &\\ F_q &&&& ~F_p  \\ &   \nwarrow & &   \nearrow & \\&& F && \\   &&  \uparrow  & \   & \\ && K && \end{array} ,
 \; \rho_{G,s}(\gamma) =  \left(\begin{array}{ccccc} 
1 & ^t\xi_q(\gamma) & \tau_s(\gamma) \\ 0 & \rho_E(\gamma) & \xi_p(\gamma) \\ 0& 0 & 1    \end{array}\right),
 \begin{array}{ccc}
              &\tau_s :    Gal_\partial(L/F_{pq}) \hookrightarrow \C \qquad&  \\
              &^t\xi_q :   Gal_\partial(F_q/F) \hookrightarrow \C^2 \simeq \hat  {\cal V}&\\
              &\xi_p :    Gal_\partial(F_p/F) \hookrightarrow \C^2 \simeq {\cal V}&\\
              &\rho_E :  Gal_\partial(F/K) \hookrightarrow SL_2(\C) & \\ 
\end{array}  $$
For convenience, we recall from \S 5.1 that the field $F = K(\omega_1, \omega_2, \eta_1, \eta_2)$ is the Picard-Vessiot extension of $K$ given by the Picard-Fuchs equation $\partial_{Lie \tilde E}(*) = 0$ for $E/S$, whose set of solutions we denote by ${\cal V} \simeq \C^2$.  If $E/S$ is isoconstant,   $F = K$,  while $Gal_\partial(F/K) = SL_2(\C)$ otherwise.  The field $F_p = F(u, \zeta(u))$ corresponds to the inhomogeneous equation attached to $p$ (given by $\partial_{Lie \tilde E}(\tilde u) = \partial \ell n_{\tilde E} (\tilde p)$, for any choice of a lift $\tilde p \in \tilde E(S)$ of $p$), while  $F_q = F(v, \zeta(v))$ is the field of periods of the semi-abelian scheme $G/S$, generated by the solutions of $\partial_{Lie \tilde G}(*) = 0$. As already said, its resemblance with $F_p$ reflects a duality, witnessed by the dual $\hat {\cal V}$ of $\cal V$. We fix a polarization of $E/S$, allowing us to identify $E$ and $\hat E$ (and in particular, $q$ to a section of $E/S$), but will keep track of this duality. The field of generalized periods of $\{E, p, q\}$ is the compositum  $F_{pq}$ of $F_p$ and $F_q$. Finally, 
$$L = F_{pq}(log_G s(\lambda)) = F_{pq}(log_{\tilde G}(\tilde s(\lambda))$$
 is the Picard-Vessiot extension generated by the solutions of the 3rd order inhomogeneous equation   $\partial_{Lie \tilde G}(\tilde U) = \partial \ell n_{\tilde G} (\tilde s)$, where $\tilde s$ is the pullback of $\tilde p$ to $\tilde G$ over 
 $s$. The corresponding 4th order homogeneous system can  be described as the Gauss-Manin connection  attached to the smooth one-motive $M$ over $S$ given by the section $s \in G(S)$. The drawing in the middle is a representation $\rho_{G,s}$ of the differential Galois group of $L/K$. The right part expresses that  the coefficients of this representation become injective group homomorphisms on the indicated subquotients of $Gal_\partial(L/K)$.  
 
 \medskip
 We will again distinguish between several cases, depending on the position of $p$ and $q$ with respect to the projection $E^\sharp$ to $E$ of the Manin kernel $\tilde E^\sharp$ of $\tilde E$. So,
\[ E^\sharp = E_{tor} + {\cal E}_0(\C) = 
       ~ \left\{ \begin{array}{ccccccc}
                			E_{tor}~ {\rm ~if}~ E ~Ê{\rm is ~not ~isoconstant ;}\\
              E_0(\C) ~{\rm ~if}~ E \simeq E_0 \times SÊ \\
                			  \end{array} \right. \\
			  \]
depending on whether the $K/\C$-trace ${\cal E}_0$ of $E$ vanishes or not. (In fact, $E^\sharp
$ is the Kolchin closure of $E_{tor}$, and is also called the {\it Manin kernel of $E$}.) We denote by $\hat p, \hat q$ the images of $p, q$ in the quotient $E/E^\sharp$. Notice that the ring ${\cal O} = End(E/S)$ still acts on this quotient.

\medskip

We recall that we must here merely prove Part (A) of the Main Lemma. By contraposition, we assume that no constant translate of $s$ is a Ribet section, or factors through a strict subgroup scheme of $G/S$, and must deduce that  $log_G(s)$, or equivalently, $log_{\tilde G}(\tilde s)$, is transcendental over $F_{pq}$.

\subsubsection*{Case (SC1) : $\hat q = 0$}

Assume first that $E/S$ is not isoconstant. Then, this vanishing means that $q$ is  a torsion section, and  after an isogeny, $G = G_q$ is isomorphic to $\G_m \times E$. We have already proven Part (A) of the Main Lemma in this case, see Lemma 1 and the lines which follow. So, we can assume that $E = E_0 \times S$ is constant, and the relation $\hat q = 0$ now means that $q$ is constant. So, $G = G_0 \times S$ is a constant semi-abelian variety, and we can apply to its (constant) universal vectorial extension $\tilde G_0$ the  slight generalization of Ax's theorem given in \cite{Be}, Proposition 1.b. Since we are assuming that no constant translate   $s - s_0, s_0 \in G_0(\C)$, of $s$  factors through a strict subgroup scheme $H$ of $G$,  the relative hull $G_s$ of $s$ in the sense of \cite{Be}, \S 1, is equal to $G$,  and  Proposition 1.b of {\it loc. cit.}  implies that 
$$tr.deg.(K(log_{\tilde G}(\tilde s))/K) =   dim(\tilde G) = 3.$$ 
Now,  $F_q = F = K$ since $E$ and $q$ are constant, while  $K(\tilde U) =K(\tilde u,  U) =  F_p(log_G(s))$ has   transcendence degre  $\leq 1$  over $F_p$, which has transcendence degree $\leq 2$ over $K$. So, both transcendence degrees must be maximal, and $log_G(s)$ is indeed transcendental over $F_p = F_{pq}$.

\subsubsection*{Case (SC2) : $\hat p = 0$}

This case is dual to the previous one, and the following preliminary remarks will simplify its study. The hypothesis made on $s$ implies that $p$ is not a torsion section (otherwise, a multiple of $s$ factors through $\G_m$). So, we can assume that $E = E_0 \times S$ is constant, and that $p = p_0$ is a constant non-torsion section of $E$. In view of  Case $({\bf CS1})$, we can also assume that $\hat q \neq 0$, i.e.  that   $q$ is not constant\footnote{~Ê
 It is worth noticing that this case $({\bf SC 2})$  is the logarithmic analogue of the counterexample  studied in \cite{BP}, \S 5.3. It does not provide a counterexample to the Main  lemma, whose ``exponential" analogue would amount, in the notations of \cite{BP}, to the    equality  $tr.deg. (K(y)/K) = 1$. In fact, the Theorem stated in \S 4, combined with Lemma 1 and with the conclusion of  Case ({\bf SC2}), implies that $y$ is transcendental over $K^\sharp_G$.}
 (and in particular, not torsion). We now consider the smooth one-motive $M = [M_s : \Z \rightarrow G]$ attached to the section $s$ above $p = p_0$, and its Cartier dual $M' = [M'_{s'} : \Z \rightarrow G']$, where $G'$ is the extension of $\hat E$ by $\G_m$ parametrized by $p_0$. In particular, $G'$ is a constant and non isotrivial semi-abelian variety. By  Remark 3.ii of \S 5,  the field  $F_{p_0q}(log_G(s)) = F_q(log_G(s))$ coincides with the field $F_{qp_0}(log_{G'}(s'))= F_q(log_{G'}(s'))$. But  the section $s'$ of $G'$ project to $q$ in $\hat E$, which is not constant, so no constant translate of $s'$  factors through $\G_m$. Finally, $\G_m$  is the unique connected subgroup scheme of $G'$, since $G'$ is non isotrivial. So, the constant semi-abelian variety $G'$ and its section $s'$ satisfy all the hypotheses of Case ({\bf SC1}). Therefore, $log_{G'}(s')$ is transcendental over $F_q = F_{qp}$, or equivalenty, $log_G(s)$ is transcendental over $F_{pq}$.

\bigskip
In the next two cases, the proof of our transcendence claim   can be derived from the following simple observation : the Lie algebra ${\mathfrak u}_s$ of the unipotent radical of the image of $\rho_{G,s}$ consists of matrices of the form $X$ indicated below, where $(^ty,x) \in Im\big((^t\xi_q, \xi_p)\big) \subset  \hat{\cal V}Ê\times {\cal V} $, and $t \in \C$, and for two such matrices
$$ X = \left(\begin{array}{ccccc} 
0& ^ty  &t  \\ 0 &0 &x \\ 0& 0 & 0    \end{array}\right), X', ~Ê{\rm we ~have}~
 [X, X'] =  \left(\begin{array}{ccccc} 
0& 0 & t(X, X')  \\ 0 &0 &0 \\ 0& 0 & 0    \end{array}\right) ,$$
where $ t(X, X') = \; <y|x'> - <y'|x>$  depends only on the vectors $x, y, x', y'$. Here, the transposition  and the scalar product represent the canonical antisymmetric pairing $ {\cal V}Ê\times {\cal V}Ê\rightarrow \C$ provided by the chosen principal polarisation on $E/S$. Now,

\subsubsection*{Case (SC3) :  $\hat p$ and $\hat q$ are linearly independent over $\cal O$}

As mentioned in Remark 2.(iv) of \S 4,  the argument  leading to Lemma 1, i.e. the sharpened form of Andr\'e's theorem \cite{An} given in \cite{BM} (or alternatively, if  $E/S$ is isoconstant, the sharpened forms of Ax's theorem given in \cite{BK}, Thm. 2,  and in \cite{Be}, Prop. 1.b)  implies in this case that $u, \zeta(u), v, \zeta(v)$ are algebraically independent over $F$. In other words, the homomorphism $(^t\xi_q, \xi_p) : Gal_\partial(F_{pq}/F) \rightarrow \hat {\cal V}Ê\times   {\cal V}Ê\simeq \C^4$ is bijective, and any  couple $(^ty,x)$ occurs in the Lie algebra ${\mathfrak u}_s$. Consequently, there exists $X, X' \in    {\mathfrak u}_s$ such that $t(X, X') \neq 0$, and ${\mathfrak u}_s$ contains matrices all of whose coefficients, {\it except} the upper right one, vanish. Therefore, the homomorphism $\tau_s $ is bijective, $Gal_\partial(L/F_{pq}) \simeq \C$, and $tr.deg(L/F_{pq}) = 1$. (Notice that 
 this yields $tr.deg (L/F) = 1 + 4 = 5$.) 
 
 \medskip
 So, from now on, we can assume that $\hat q$ and $\hat p$ are linked by a unique relation over $\cal O$, which, considering multiples if necessary, we write in the shape
 $$\hat q = \alpha \hat p, \alpha \in {\cal O}, \hat p \neq 0, \alpha \neq 0.$$
 We denote by $\overline \alpha$ the complex conjugate of $\alpha
 $, which represents the image of $\alpha \in End(E/S)$ under the Rosati involution attached to the chosen polarization. We first deal with  non antisymmetric relations, in the sense of \S 2.3.
 
 \subsubsection*{Case (SC4) :  $\hat q = \alpha \hat p$, where $  \overline \alpha  \neq - \alpha$ and $\hat p \neq 0$}  
 
Lifted to $E$ and up to an isogeny, this relation reads  $q = \alpha p + p_0$, where $p_0 \in E_0(\C)$ is a constant section, equal to $0$ is $E/S$ is not isoconstant.    Then, $F_q = F_p$, and more precisely, $\xi_q = \alpha \circ \xi_p$. In other words, the coefficients of the matrices $X$ in ${\mathfrak u}_s$ satisfy the relation $y = \alpha x$, for the natural action of $\cal O$ on $\cal V$.  Since $p \notin E_0(\C)$, Lemma 1 implies that $\xi_p$ is bijective, and any $x \in {\cal V}$ occurs  in the Lie algebra ${\mathfrak u}_s$.  Recall that $< | >$ is antisymmetric, and that the adjoint of the endomorphism of $\cal V$ induced by an isogeny $\alpha \in {\cal O}$ is its Rosati image $\overline \alpha$. For any $x, x'$ in $\cal V$ occurring in matrices $X, X'$ and such that $<x|x'> \,\neq \, 0$,  we then have :
 \begin{eqnarray*} t(X, X') &=&   < \alpha x | x'> - <\alpha x' |x> \; =\;  < \alpha x | x'> - <  x' | \overline \alpha x> \\
&=&  < \alpha x | x'>  + <\overline \alpha x, x'> \; = \;  <(\alpha + \overline \alpha)x | x'> \; \neq \; 0   
\end{eqnarray*} 
since $\alpha + \overline \alpha$ is a non-zero integer. We conclude as in Case $({\bf SC3})$ that $Gal_\partial (L/F_{pq}) \simeq \C$, and $tr.deg.(L/F_{pq}) = 1$. (Here, this yields $tr.deg(L/F) = 1 + 2 = 3$.)

\medskip
The remaining cases concerns antisymmetric relations of the type $\hat q = \alpha \hat p$, with $\hat p \neq 0$ and a non zero purely imaginary $  \alpha = - \overline \alpha$.  (In particular, the CM elliptic scheme $E/S$ must be isoconstant and so, Theorem 2 is now already proven under its Conditions (i) or (ii)). We first treat the case when $q$ and $p$ themselves are antisymmetrically related.

 \subsubsection*{Case (SC5) :  $ q = \alpha  p$, where $  \overline \alpha  = - \alpha \neq 0$ and $\hat p \neq 0$} 
 
 This is the only case where a Ribet section of $G/S$ exists above the section $p \in E(S), p \notin E_0(\C)$. Denote by $s_R$ this (essentially unique) Ribet section. 
 
\begin{Lemma} Let $s_R$ be the Ribet section of $G/S$. Then, $log_G(s_R)$ is defined over $F_{pq}$.
\end{Lemma}
\noindent
{\it Proof} : let $L_R = F_{pq}(log_G(s_R))$ be the field generated over $F_{pq}$ by $log_G(s_R)$. Since the differential Galois group $Gal_\partial(L_R/F_{pq})$ injects via $\tau_{s_R}$ into a vectorial group $\C$, $log_G(s_R)$ is either transcendental or rational over $F_{pq}$. Assume that it is transcendental.  Then, by Proposition 4, the surface $\cal S$ attached to $s_R$ contains no algebraic curve, and the whole reduction of the Main Theorem to the study of $\cal S$  given in \S 3.1 and \S 3.2   implies that $s_R$ admits only finitely many torsion values. But this contradicts the main result of \cite{BE}, an analytic proof (or version) of which will be found in   Appendix I, \S 7.2. In fact, the explicit formulae given there directly show that $L_R = F_{pq}$, cf. Propositiion 5.(iii).

\medskip
We now come back to our section $s$, which we assumed {\it not} to be a Ribet section of $G$, and more accurately, such that no constant translate $s- s_0$ of $s$ is a Ribet section. Since $s$ and $s_R$ project to the same section $p$ of $E/S$, there exists a section $\delta \in \G_m(S)$, i.e. a rational function in $K^*$, such that $s = s_R + \delta$. The assumption on $s$ implies that $\delta \notin \C^*$ is not constant. Set $\ell = log_{\G_m}(\delta)$. Then, $log_G(s) = log_G(s_R) + \ell$, and the lemma implies that that $F_{pq}(log_G(s)) = F_{pq}(\ell)$, which is equal to $F_p(\ell)$, since $\hat q = \alpha \hat p$. By Lemma 1, $F_p(\ell)$ has transcendence degree 1 over $F_p$, so $log_G(s)$ too is transcendental over $F_{pq}$.

\subsubsection*{Case (SC6) :  $ q = \alpha  p + p_0$, where $  \overline \alpha  = - \alpha \neq 0$ and $p_0 \in E_0(\C), p_0 \notin E_{0, tor}$}

 Fixing an $\alpha$-division point $p'_0$ of $p_0$ in $E_0(\C)$, we have $ q = \alpha (p- p'_0)$, and there exists an (essentially unique) Ribet section $s'_R$ of $G/S$  above $p' := p -p'_0$. Then, the section $s' := s -s'_R$ of $G/S$ projects to $\pi(s') = p - p' = p'_0$ in $E_0(\C) \subset E(S)$. Furthermore, $F_p = F_{p'} = F_q$ since $p'_0$ is constant, and  since $log_G(s'_R)$ is defined over $F_{p'q} = F_{pq} = F_q$ by Lemma 3, we deduce that $log_G(s') = log_G(s)- log_G(s'_R)$ generates over the field $F_{pq} = F_{p'q} = F_q = F_{p'_0q}$ the same field as $log_G(s)$.  We are therefore reduced to showing that given $\hat q \neq 0$ and a section $s'$ of $G/S$ projecting to a constant non-torsion section $p'_0 \in E_0(\C)$, then $log_G(s')$ is transcendental over $F_{p'_0q}$. But this is exactly what we proved in Case $({\bf SC2})$ !

 \smallskip
This concludes the proof of the Main Lemma, hence of the Main Theorem.

 \section{Appendix I}

\subsection{Analytic description of the semi-abelian logarithm}

Let $G/S$ be a {\it non isotrivial} extension of an elliptic scheme $E/S$ by $\G_m$ and let $s$ be a section of $G/S$.  The aim of this Appendix is to give an explicit formula for its local logarithm $log_G(s)$ in terms of the  Weierstrass functions $\wp_\lambda, \zeta_\lambda, \sigma_\lambda$, in parallel with that of \S 4 in for products. We recall  the notations of \S 5.1. In particular,   we set $p = \pi(s) \in E(S)$, $u = log_E(p)$, For simplicity, we will work over the generic point of $S$,  consider $G$ as a semi-abelian variety over the field $K = \C(S)$, and  often drop the variable $\lambda$  indexing the  Weierstrass functions and their (quasi-)periods  $\omega_1, \omega_2, \eta_1, \eta_2$.

\medskip
By Weil-Rosenlicht-Barsotti, the algebraic group $G$, viewed as a $\G_m$-torsor, defines a line bundle over $E$ of degree 0, admitting a rational section  $\beta$ with divisor $(-q) - (0) \in \hat E$, which we identify with the point $q \in E$ (the sign is admittedly not standard, but it will make the formulae symmetric in $p$ and $q$).  By assumption, $q$ is not a torsion point, and we set $v = log_E(q)$. We further assume that $p \neq 0$ and $p + q \neq 0$. 

\medskip

The rational section $\beta$ provides a birational   isomorphism  $G \dasharrow  \G_m \times E$ and (after a shift away from $0$) an isomorphism $Lie G \simeq Lie\G_m \oplus Lie E$. The 2-cocycle which describes   the group law on the product (\cite{Se}, VII.5) is a rational function on $E \times E$, expressed in terms of $\sigma$-functions by $\frac{\sigma(  z + z' + v) \sigma(z) \sigma(z') \sigma(v)}{\sigma(z+z') \sigma(z+v) \sigma(z'+v)}$. Therefore, the exponential morphism $exp_G$ is represented by the map
$$(Lie G)^{an}(\Lambda) \ni  \left( \begin{array}{c} t(\lambda) \\z(\lambda) \end{array}   \right)  \mapsto  \left( \begin{array}{c} f_{v(\lambda)}(z(\lambda))\;e^{t(\lambda)} \\ \wp_\lambda(z(\lambda)) \end{array}   \right)   \in G^{an}(\Lambda)  $$
where  
    $$f_v(z) = \frac{\sigma (v+z)}{\sigma (v) \sigma (z)} e^{-\zeta(v)z} $$
 is a meromorphic  theta function for the line bundle  ${\cal O}_{E}\big((-q)- (0)\big)$, whose factors of automorphy are given by $e^{-\kappa_v(\omega_i)}$ (opposite of the multiplicative quasi-periods), with 
$$\kappa_v(\omega_i) =  \zeta (v)   \omega_i - \eta_i v , ~Ê{\rm for}~ i = 1,2.$$ 
The occurence of the trivial theta function $e^{-\zeta(v) z}$ in $f_v$ is due to the condition $d_0(exp_G) = id_{Lie G}$. The logarithmic form $\frac{df_v}{f_v} = \big(\zeta(v+z) - \zeta(v) - \zeta(z)\big) dz  = \frac{1}{2} \frac{\wp'(z) - \wp'(v)}{\wp(z) - \wp(v)} dz$ is the pullback under $exp_E$ of the standard differential form  of the 3rd kind on $E$ with residue divisor $-1.(0) + 1.(-q)$. 

\medskip
 Under this description, the   section $s$ of $G/S$ under study and its logarithm $log_G(s) $ are  given by
$$s =  \left( \begin{array}{c} 
 \delta_s    \\ p  \end{array}   \right) , ~ U  := log_G(s) =  \left( \begin{array}{c} 
-g(u,v) + \zeta(v)u   + \ell_s    \\ u  \end{array}   \right)$$
where   $\delta_s := s - \beta(p) \in K^*$ is   a rational function on $S$, depending only on $s$ (and on the choice of the section $\beta$), for which we set $\ell_s = log_{\G_m}(\delta_s)$,  and (cf. the formulae in \cite{Be}, up to signs) : 
$$ g_\lambda(u,v) = log\big(\frac{\sigma_\lambda(u+v)}{\sigma_\lambda(v) \sigma_\lambda(u)}\big).$$
This is the Green function mentioned in \S 5.1.

\medskip

The $\Z$-local system  of periods $\Pi_G$ of $G^{an}/\Lambda$ which was introduced in \S 3.3  admits the basis
$$\varpi_0(\lambda) =   \left( \begin{array}{c} 
  2 \pi i \\ 0  \end{array}   \right) , \varpi_1(\lambda) =   \left( \begin{array}{c}   \kappa_{v(\lambda)}(\omega_1(\lambda)) \\   \omega_1(\lambda)  \end{array}   \right) , \varpi_2(\lambda) =  \left( \begin{array}{c}  \kappa_{v(\lambda)}(\omega_2(\lambda)) \\   \omega_2 (\lambda)   \end{array}   \right). $$
  
  \medskip
  
  We can now describe the various extensions of $F = K(\omega_1, \omega_2, \eta_1, \eta_2)$  appearing in \S 5 for a non isotrivial extension $G$.  In view of the Legendre relation
  $$2 \pi i = \eta_1 \omega_2 - \eta_2 \omega_1 \in K^*, $$
  the  periods of $G$ generate over $F$ the field  
  $$F_G  = F^{(2)}_G  = F(\kappa_v(\omega_1), \kappa_v(\omega_2)) =   F( v, \zeta(v)) := F_q,$$
 while the field generated over $F_{pq}$ by $log_G(s)$ satisfies
 $$L = F_{pq}(log_G(s)) = F(u, \zeta(u), v, \zeta(v), -g(u,v) + \zeta(v) u+ \ell_s) = F_{pq}\big(\ell_s -g(u,v)\big)$$

\subsection{Analytic description of the Ribet sections}
\medskip
 
We now  present the  analytic  description of the   Ribet sections promised in \S 1.1.  Assume that  $E = E_0 \times S$ is a {\it constant elliptic scheme with complex multiplications}, that {\it $q$ is not constant} (i.e. $G$ is semi-constant), and that  {\it $p$ and $q$ are antisymmetrically related} in the sense of \S 2. So, their logarithms   $u(\lambda), v(\lambda)$ are non constant holomorphic functions on $\Lambda$, and satisfy~:  $v = \alpha u$ modulo $\Q \omega_1 \oplus \Q \omega_2$ for a totally imaginary non zero complex multiplication $\alpha \in {\cal O}Ê\otimes \Q$. For simplicity, we will assume that this relation takes the form
$$ v = \alpha u , \alpha \in {\cal O}, \alpha  = - \overline \alpha \neq 0, u, v \notin \C,$$
and that $\alpha$ lies in $2 {\cal O} \subset 2( \Z \oplus \Z \tau)$, where $\tau := \frac{\omega_2}{\omega_1}Ê\in {\cal O}Ê\otimes \Q$, $\Pi_E =  \omega_1(\Z \oplus \Z \tau)$.

\medskip

We fix a rational section $\beta$ of $\pi : G \rightarrow E$ as above, so, any section $s$ of $G/S$ projecting to $p \in E(S)$ is uniquely expressed as a couple $(\delta_s, p)$, for some rational function $\delta_s \in  K^*$. In particular, the Ribet section $s_R$ of $G/S$ is given by $(\delta_R, p)$, where $\delta_R \in K^*$ is well defined, at least up to a root of unity. Since $E$ is constant,  the index $\lambda$ truly disappears from the  notations $\omega_{1,2}, \sigma, g, ...$, which are now attached to $E_0$.  We must express $\delta_R$ in terms of these functions, evaluated at $u(\lambda)$.

\medskip
By the theory of complex multiplication  (see for instance  \cite{Ma}, Appendix I), there exists  an explicit complex number $s_2$ (actually algebraic over $\Q$, and given by  Hecke's non-holomorphic modular form  of weight 2), such that the quasi-periods of $\zeta$ satisfy :
$$   \eta_2 - s_2 \omega_2 = \overline \tau (\eta_1 - s_2 \omega_1), ~Ê{\rm where}~ \omega_2 = \tau \omega_1.  $$
From this and the functional equation of $\sigma$, one infers that   the function $\theta(z) = \sigma(z) e^{- \frac{1}{2}s_2z^2}$ satisfies,  for  $N(\gamma) = \gamma \overline \gamma$ :
$$ \forall \gamma \in {\cal O}~, ~ \Big(\frac{\theta(\gamma z)}{\theta(z)^{N(\gamma)}}\Big)^2 =    \gamma^2   \prod_{e \in E[\gamma], e \neq 0} (\wp(z) - \wp(e)).$$
Since $N(\alpha + 1) - N(\alpha) - 1 = \alpha + \overline \alpha = 0$, we deduce from the parity assumption on $\alpha$ that the  meromorphic function
$$\delta_{\tilde R} := e^{g(u,  v )  - s_2u    v} = \frac{\sigma(u +  v )}{\sigma ( v ) \sigma(u)} e^{-s_2 u   v},  ~Ê{\rm where}Ê~ v = \alpha u,$$
 is a rational function of $\wp(u(\lambda)), \wp'(u(\lambda))$, hence  lies in $K^*$. We then have:
 \begin{Proposition} For $q = \alpha p$, consider the section $s_{\tilde R} = (\delta_{\tilde R}, p)$  of $G = G_q$ above $p$.  Then,
 
 i) up to a root of unity, $s_{\tilde R}$ is equal to the Ribet section $s_R$ of $G$ above $p$;
 
 ii) for any $\overline \lambda \in S$ such that $p(\overline \lambda)$ is a torsion point of $E_0$, say of order $n$, the point $s_{\tilde R}(\overline \lambda)$ is a torsion point of $G_{\overline \lambda}$, of order dividing $n^2$;
 
 iii) the logarithm of $s_{\tilde R}$ is given by $(\zeta(v) u - s_2uv, u)$; in particular, $F_{pq}(log_G(s_R)) = F_{pq}$.
 
 \end{Proposition}
\noindent
(Following the conventions of \S 2.1, the bar over a $\lambda$ indicates a value of $\lambda$, not complex conjugation as for $\tau, \alpha$ and $\gamma$.)  

\bigskip
\noindent
{\it Proof}. -  i) In view of Corollary 3, applied to $s = s_{\tilde R}$ and $s^\dagger = s_R$, this first assertion is an immediate corollary of the second statement (ii). Indeed, by Proposition 1.(iv), the set $S_{\infty}^E$ is infinite since $p$ is not constant, so (ii) states in particular that  $s_{\tilde R}$ lifts infinitely many torsion values of $p$ to torsion points of $G$.

\medskip

iii) By the formulae of the previous subsection, the first coordinate of $log_G(s_{\tilde R})$ is equal to
$$- g(u,v) + \zeta(v) u + log(\delta_{\tilde R}) = \zeta(v) u -s_2 uv.$$
The last assertion then immediately follows from (i). It is timely to point out here that the logarithmic derivative $\tilde \zeta(z) = \zeta(z) - s_2 z$ of $\theta$  satisfies : $\tilde \zeta(\alpha z) - \overline \alpha \,\tilde \zeta(z) \in \wp' (z)\C( \wp(z))$.

\medskip

ii) We must  show that for  any $\overline \lambda$ such that  $u(\overline \lambda) = \frac{1}{n} \omega$ for some integer $n$ and some period $\omega$ of $E$, with corresponding additive and multiplicative quasi-periods $\eta(\omega),   \kappa_v(\omega)$, and with $v(\overline \lambda) =  \alpha u(\overline \lambda) = \frac{1 }{n}\alpha \omega$,  the number
$$\zeta(v(\overline \lambda)) u(\overline \lambda) - s_2 u(\overline \lambda)v(\overline \lambda) -  \frac{1}{n} \kappa_{v(\overline \lambda)}(\omega) = \frac{1}{n^2} \eta(\omega) \alpha \omega    -  \frac{1}{n^2} s_2 \alpha  \omega^2= \frac{1}{n^2} \alpha (\eta(\omega) - s_2 \omega) \omega $$
is a rational multiple of  $2\pi i$, and more precisely an integral multiple of $\frac{1}{n^2}.2\pi i$. Indeed, this number is the first coordinate of  the vector $log_G(s_{\tilde R}(\overline \lambda)) - \frac{1}{n} \varpi$, where $\varpi$ lies in the $\Z$-module $\Pi_G$ of periods of $G$ above the period $\omega \in \Pi_E$. Since $\log_E(p(\overline \lambda)) = u(\overline \lambda) =  \frac{1 }{n} \omega$, the second coordinate of this vector vanishes, so  $log_G(s_{\tilde R}(\overline \lambda))$  itself will lie $\frac{1}{n^2} \ZÊ\varpi_0 + \frac{1}{n}Ê\Pi_G$.  

\smallskip

By additivity, it suffices to  check the assertion when $\omega = \omega_1$ and  $\omega_2 = \tau \omega_1$. By the classical Legendre relation, combined with the CM relation $\eta_2 - s_2 \omega_2 = \overline{\tau} (\eta_1 - s_2 \omega_1)$,  the number under study is equal to $\frac{1}{n^2} \frac{\alpha}{\overline \tau - \tau} 2\pi i$ when $\omega = \omega_1$. Since $\alpha$ is a purely imaginary element of $\Q(\tau)$, this is indeed a rational multiple of $2 \pi i$, and even an integral one since $\alpha \in \Z \oplus 2\Z \tau$.  For $\omega = \omega_2$, the same computation yields $\frac{1}{n^2} \tau \overline \tau \frac{\alpha}{\overline \tau - \tau} 2\pi i$, with same conclusion since $\alpha \tau \in 2\Z \oplus \Z \tau$. 

\bigskip
\noindent
{\bf Remark 4.}  i) Proposition 5.(ii) by itself suffices to give a counterexample to the Relative Manin-Mumford conjecture, independently of \cite{BE}. So, we could have  used $s_{\tilde R}$ as a definition of the Ribet section. Granted the equality $s_{\tilde R} = s_R$,  this gives a new proof of Edixhoven's sharpening in \cite{BE}, Appendix, according to which the Ribet section $s_R$ lifts torsion points of order $n$ to torsion points of order dividing, and often equal  to,  $n^2$.   

\smallskip
ii) Instead of a parity assumption, one can merely assume that $\alpha = \beta - \overline \beta$ for some $\beta \in \cal O$. As shown in greater generality in \cite{B-E}, the $n$-th root of unity   $n s_R(\overline \lambda)$, for $\overline \lambda \in S^E_n$, can then be expressed in terms of the Weil pairing $e_n(\beta p (\overline \lambda), p (\overline \lambda))$. The formula above also implies this sharpening, since it shows that $n s_{\tilde R}(\overline \lambda) $ is then equal to $exp(\frac{1}{n} \big( \eta(\omega) \omega' - \eta(\omega') \omega) \big)$, where $\omega' = \beta \omega$. 

\section{Appendix II : application to Pell equations}

 In \cite{MZbis}, the relative Manin-Mumford conjecture is proven for simple abelian surface schemes, and this implies the following corollary :  consider the family of sextic polynomials  $D_\lambda(x) = x^6 + x + \lambda$, where $\lambda$ is a complex parameter. Then, there are only finitely many  $\overline \lambda \in \C$ such that the functional Pell equation $X^2 - D_{\overline \lambda}(x) Y^2 = 1$ admits a solution in   polynomials $X , Y \in \C[x]$,  $Y \neq  0$; see also \cite{Za}, III. 4.5 for connections with other problems and a proof of the deduction from RMM . The involved abelian surface $A/\C(\lambda)$  is the jacobian of the (normalized) relative hyperelliptic curve $C : y^2 = x^6 + x + \lambda$, and RMM is applied to the section $s$ of $A$ defined by the linear equivalence class of the relative divisor $(\infty_+) - (\infty_-)$ on   $C$.
 
 \medskip
 Following a suggestion of the second author, we may treat in the same way the case of a sextic $D_\lambda(x) = (x- \rho(\lambda))^2Q_\lambda(x)$   having a {\it squared linear factor}, i.e. a generic double root $\rho(\lambda)$ for some algebraic function $\rho(\lambda)$, now applying the Main Theorem of the present paper to a quotient $G= G_\rho$ of the generalized jacobian  of the corresponding semi-stable relative sextic curve $C$. This $G_\rho$ is an extension by $\G_m$ of an elliptic curve  $E$ (over $\C(\lambda)$), where $E$ is the jacobian of the (normalized) relative quartic $\tilde C$   with equation $v^2 = Q_\lambda(u)$, and RMM may  be applied to the section $s$ of $G_\rho$ defined by the  class  of the relative divisor $(\infty_+) - (\infty_-)$ on  $\tilde C$, for the strict linear equivalence attached to the node of $C$ at $x = \rho(\lambda)$ (see \cite{Se}, V.2, \cite{BE}, Appendix). When the quartic $Q_\lambda(x)$ does not depend on $\lambda$ (or more generally, when $ E$ is isoconstant), no result of interest seems to appear. As an illustration of the opposite case, we will here consider the family of quartics
 $$Q_\lambda(x) = x^4 + x + \lambda.$$
From the analysis in \cite{MZbis}, \cite{Za} recalled below, we derive that the set $\Lambda_Q$ of complex numbers $\overline \lambda$ such that the Pell equation $X^2 - Q_{\overline \lambda}(x) Y^2 = 1$ has a solution  in polynomials  $X , Y \in \C[x]$,  $Y \neq 0$, is infinite. The solutions for each such $\overline \lambda$ form a sequence $(X_{\overline \lambda, n}, Y_{\overline \lambda, n})_{n \in \Z}$ of polynomials in $\C[x]$. Our result is that for the $\rho(\lambda)$ considered in the following statements,  only finitely many of  the polynomials $Y_{\overline \lambda, n}(x)_{\overline \lambda \in \Lambda_Q, n \in \Z}$     admit $x = \rho(\overline \lambda)$ among their roots. In other words : 
 
 \begin{Theorem} i) (Case $\rho(\lambda) = 0$.) There are only finitely many complex numbers $\overline \lambda$ such that the equation $X^2 - x^2Q_{\overline \lambda}(x) Y^2 = 1$ admits a solution in polynomials $X, Y \in \C[x], Y \neq 0$.
 
 ii) (Case $\rho(\lambda) =  \frac{4y(\frac{1}{2}p_W + e_3) - 1}{8x(\frac{1}{2}p_W + e_3) } \in \C(\lambda)^{alg} $, with notations explained below.)  There are only finitely many complex numbers $\overline \lambda$ such that the equation $X^2 - (x- \rho(\overline \lambda))^2Q_{\overline \lambda}(x) Y^2 = 1$ admits a solution in polynomials $X, Y \in \C[x], Y \neq 0$.

 \end{Theorem}
 
 In spite of their similarity, these two statements cover different situations : in (i), the extension $G_\rho$ is not isotrivial, and the theorem is a corollary of Theorem 2. On the contrary, (ii) illustrates the case of an isotrivial extension $G_\rho$, and follows from Theorem 1. In fact, we believe that on combining these two cases of our Main Theorem, Theorem  4 will hold for {\it any} choice of the algebraic function $\rho(\lambda)$. 
 
 \medskip
 
The algebraic function $\rho(\lambda)$ of Case (ii) can be described as follows. Consider the Weierstrass model
 $(W_E) : y^2 = 4 x^3 - \lambda x + (1/16)$ of the elliptic curve $E/\C(\lambda)$, with its standard group law, and the relative point $p_W =  (0, -\frac{1}{4})$ on $(W_E)$, which is generically of infinite order (it can be rewritten as the point $p_{\tilde W} = (0, -1)$ on the curve   $(\tilde W_E) : Y^2 = X^3 - 4 \lambda X + 1$, and is not torsion at $\overline \lambda = \frac{1}{4}$). Then, the $2$-division points of $p_W$ are  the four points of $(W_E)$ :
 $$`` \frac{1}{2}p_W " = (\frac{1}{8}m_\lambda^2, -\frac{1}{8}m_\lambda^3 + \frac{1}{4}),~Ê{\rm  ~ where} ~ m_\lambda  ~Ê{\rm ~ is ~a ~root ~of} ~ m^4 - 8m + 16 \lambda = 0.$$
 Choose one of the two roots $m_\lambda$ which is real when $\overline \lambda = \frac{1}{4}$, and call the corresponding point $\frac{1}{2} p_W(\lambda)$. Further, choose one of the two points of order 3 on $W_E$ which is real when $\overline \lambda$ is real, and call it $e_3(\lambda)$. Computing the $x$ and $y$ coordinates of the relative point $\frac{1}{2}p_W + e_3$ on $W_E$ then provides the function $\rho(\lambda)$ appearing in Case (ii).

 \medskip
In a more enlightening way, let in general $p(\lambda)$ be the section of $E$ defined by the class of the divisor $(\infty_+) - (\infty_-)$ on $\tilde C$, for the standard linear equivalence of divisors. Then, $p$ is the projection to $E$ of the section $s$ of $G_\rho$ defined (via strict equivalence) above, and one checks that $p$ is not a torsion section. By \cite {MZbis}, \cite{Za} (see also \cite{Haz}, Prop. 3.1), the Pell equation for $Q_{\overline \lambda}(x)$ has a non trivial solution  if and only if   $p(\overline \lambda)$ is a torsion point on $E_{\overline \lambda}$, i.e. $\overline \lambda \in S_\infty^E$ in the notations of \S 2.  Similarly, the Pell equation for $(x- \rho(\overline \lambda))^2 Q_{\overline \lambda}(x)$ has a non trivial solution   if and only if   $s(\overline \lambda)$ is a torsion point of $G_{\rho(\overline \lambda)}$,  i.e. $\overline \lambda \in S_\infty^{G_\rho}$. Furthermore, the section  $q(\lambda)$ of $\hat E$ parametrizing the  extension $G_\rho$ is represented the (standard) equivalence class of the divisor $(q_+) - (q_-)$ on $\tilde C$, where $q_\pm(\lambda)$ is the section $(\rho(\lambda), \pm Q^{1/2}_\lambda(\rho(\lambda)) $ of $\tilde C$. Now, in the first case $\rho(\lambda) = 0$, $q$ is a non torsion section, i.e. $G_\rho$ is a non isotrivial extension of the  non isoconstant elliptic scheme $E$, and  since $p$ is not torsion, Theorem 2 implies that $S_\infty^{G_\rho}$ is finite. On the other hand,   Case (ii) is built up in such a way that $q$ has finite order (equal to 3), so that $G_\rho$ is now isogenous to $\G_m \times E$. But one can check (by specializing at the real number $\overline \lambda = \frac{1}{4}$) that the projection of the section $s$ to the $\G_m$ factor is not a root of unity, so $s$ does not factor through a translate of $E$.  Since $p$ is not torsion either, Theorem 1 now provides the finiteness of  $S_\infty^{G_\rho}$.

 \bigskip
 \noindent
 {\it Adresses of authors: }  D. B. (Paris), D. M. (Basel), A. P. (Leeds), U. Z. (Pisa).

\end{document}